\documentclass[11pt, a4paper]{article}

\usepackage{amsmath}
\usepackage{amsfonts}
\usepackage{amssymb}
\usepackage[francais,english]{babel}

\def\english{\selectlanguage{english}}

\providecommand\mathbb{\bf}
\newcommand\R{{\mathbb R}}

\newtheorem{thm}{Theorem}[section]
\newtheorem{lemma}{Lemma}[section]

\newtheorem{pro}{Proposition}[section]
\newtheorem{defi}{Definition}[section]

\newtheorem{remark}{Remark}[section]

\newcounter{Remark}

\renewcommand\theRemark{\arabic{Remark}}

\newcounter{steps}
\newenvironment{proof}[1][]{%
\par\medbreak\setcounter{steps}{0}
{\noindent\bfseries Proof#1. }} {\hfill\fbox{\ }\medbreak}

\newcounter{substeps}[steps]

\renewcommand{\P}{
{\cal P}}

\newcommand{\Po}{
\P_1}            

\def\lip{\mathrm{Lip}}

\newcommand{\tvarphi}[0]{
\tilde{\varphi}}

\newcommand{\tpsi}[0]{
\tilde{\psi}}

\newcommand{\sphere}[0]{
\mathbb{S}}

\newcommand{\supp}[0]{
\mathrm{supp\;}}

\newcommand{\Xe}[0]{
X ^\varepsilon }

\newcommand{\Ve}[0]{
V ^\varepsilon }

\newcommand{\abv}[0]{
(\alpha - \beta |v|^2) v}

\newcommand{\abvs}[0]{
(\alpha - \beta \left |{\cal V}(s;v)\right |^2) {\cal V}(s;v)}

\newcommand{\abves}[0]{
(\alpha - \beta \left |\Ve(s)\right |^2) \Ve (s)}

\newcommand{\intxv}[1]{
\int _{\R^d \times \R^d}\!\!\! #1}

\newcommand{\intxo}[1]{
\int _{\R^d \times r \sphere}\!\!\! #1}

\newcommand{\eps}[0]{
\varepsilon}

\newcommand{\teps}[0]{
\overline{\varepsilon}}

\newcommand{\fe}[0]{
f ^\varepsilon}

\newcommand{\fek}[0]{
f ^{\varepsilon _k}}

\newcommand{\lae}[0]{
\lambda ^\varepsilon}

\newcommand{\fin}[0]{
f ^{\mathrm{in}}}

\newcommand{\fo}[0]{
f ^{(1)}}

\newcommand{\A}[0]{
\{0\} \cup r \sphere}

\newcommand{\xA}[0]{
\R ^d \times ( \{0\} \cup r \sphere )}

\newcommand{\Divx}[0]{
\mathrm{div}_x}

\newcommand{\Divv}[0]{
\mathrm{div}_v}

\newcommand{\Divo}[0]{
\mathrm{div}_\omega}

\newcommand{\ave}[1]{
\left \langle #1 \right \rangle }

\newcommand{\litwoix}[0]{
L^\infty ( \R_+ ; W^{1,\infty}(\R ^d ))}

\newcommand{\mbxv}[0]{
{\cal M}_b ^+ (\R ^d \times \R ^d)}

\newcommand{\poxv}[0]{
{\cal P}_1  (\R ^d \times \R ^d)}

\newcommand{\litmbxv}[0]{
L^\infty (\R_+ ; {\cal M}_b (\R^d \times \R ^d))}

\newcommand{\cztmbxv}[0]{
C (\R_+ ; {\cal M}_b (\R^d \times \R ^d))}

\newcommand{\cztpoxv}[0]{
C (\R_+ ; {\cal P}_1 (\R^d \times \R ^d))}

\newcommand{\imvv}[0]{
\left ( I - \frac{v \otimes v}{|v|^2} \right ) }

\newcommand{\imoo}[0]{
\left ( I - \frac1{r^2}({\omega \otimes \omega}) \right ) }

\newcommand{\czcxv}[0]{
C^0_c (\R^d \times \R^d)}

\newcommand{\czcxo}[0]{
C^0 _c (\R^d \times r \sphere)}

\newcommand{\cocxo}[0]{
C^1 _c (\R^d \times r \sphere)}

\newcommand{\czc}[0]{
C^0_c (\R^d)}

\newcommand{\cocxv}[0]{
C^1_c (\R^d \times \R^d)}

\newcommand{\ctcxv}[0]{
C^2_c (\R^d \times \R^d)}

\newcommand{\cocv}[0]{
C^1_c (\R^d)}

\newcommand{\vsv}[0]{
\frac{v}{|v|}}

\newcommand{\coctxv}[0]{
C^1_c (\R_+ \times \R^d \times \R^d)}

\newcommand{\lotczcxv}[0]{
L^1(\R_+; C^0 _c (\R^d \times \R^d))}

\newcommand{\inttxv}[1]{
\int _{\R_+} \!\int _{\R ^d \times
\R^d}\!\!\!\!\!\!#1\;\mathrm{d}t}

\newcommand{\dxv}[0]{
\mathrm{d}(x,v)}

\newcommand{\linf}[0]{
L^\infty}

\newcommand{\lime}[0]{
\lim _{\varepsilon \searrow 0}}

\newcommand{\limd}[0]{
\lim _{\delta \searrow 0}}

\newcommand{\limk}[0]{
\lim _{k \to +\infty }}

\newcommand{\reo}[0]{
\rho ^\varepsilon _1}

\newcommand{\ret}[0]{
\rho ^\varepsilon _2}

\newcommand{\reth}[0]{
\rho ^\varepsilon _3}

\newcommand{\ot}[0]{
\overline{t}}

\newcommand{\dv}[0]{
\mathrm{d}v}

\newcommand{\intv}[1]{
\int _{\R ^d} \!#1 }

\newcommand{\intvN}[1]{
\int _{\R ^d} \!#1 \;\mathrm{d}v}

\newcommand{\ind}[1]{
{\bf 1}_{#1}}

\newcommand{\psixv}[0]{
\psi \left ( x, r \frac{v}{|v|}\right )}

\newcommand{\chivd}[0]{
\chi \left ( \frac{|v|}{\delta}\right)}

\newcommand{\chipvd}[0]{
\chi ^{\;\prime}\left ( \frac{|v|}{\delta}\right)}

\newcommand{\mbxo}[0]{
{\cal M}_b ^+ (\R^d \times r \sphere)}

\newcommand{\aeps}[0]{
a ^\varepsilon}

\newcommand{\aek}[0]{
a ^{\varepsilon _k}}

\newcommand{\xp}[0]{
x^{\prime}}

\newcommand{\xs}[0]{
x^{\prime \prime}}

\newcommand{\vp}[0]{
v^{\prime}}

\newcommand{\vs}[0]{
v^{\prime \prime}}

\newcommand{\md}[0]{
\mathrm{d}}

\newcommand{\dxpvp}[0]{
\mathrm{d}(x^\prime,v^\prime)}

\baselinestretch\renewcommand{\baselinestretch}{1.5}

\begin{document}
\english

\title{Asymptotic Fixed-Speed Reduced Dynamics for Kinetic Equations in Swarming}

\author{
Mihai Bostan
\thanks{Laboratoire d'Analyse, Topologie, Probabilit\'es LATP, Centre de
Math\'ematiques et Informatique CMI, UMR CNRS 7353, 39 rue Fr\'ed\'eric Joliot Curie, 13453 Marseille  Cedex 13
France. E-mail : {\tt bostan@cmi.univ-mrs.fr}}
, J. A. Carrillo
\thanks{ICREA (Instituci\'o Catalana de Recerca i Estudis Avan\c{c}ats)
and Departament de Matem\`atiques, Universitat Aut\`onoma de
Barcelona, 08193 Bellaterra Spain. E-mail : {\tt
carrillo@mat.uab.es}. {\it On leave from:} Department of
Mathematics, Imperial College London, London SW7 2AZ, UK.}}

\date{ (\today)}

\maketitle

\begin{abstract}
We perform an asymptotic analysis of general particle systems
arising in collective behavior in the limit of large
self-propulsion and friction forces. These asymptotics impose a
fixed speed in the limit, and thus a reduction of the dynamics to
a sphere in the velocity variables. The limit models are obtained
by averaging with respect to the fast dynamics. We can include all
typical effects in the applications: short-range repulsion,
long-range attraction, and alignment. For instance, we can
rigorously show that the Cucker-Smale model is reduced to the
Vicsek model without noise in this asymptotic limit. Finally, a
formal expansion based on the reduced dynamics allows us to treat
the case of diffusion. This technique follows closely the
gyroaverage method used when studying the magnetic confinement of
charged particles. The main new mathematical difficulty is to deal
with measure solutions in this expansion procedure.
\end{abstract}



\newpage

\section{Introduction}
\label{Intro}
\indent

This paper is devoted to continuum models for the dynamics of
systems involving living organisms such as flocks of birds, school
of fish, swarms of insects, myxobacteria... The individuals of
these groups are able to organize in the absence of a leader, even
when starting from disordered configurations \cite{ParEde99}.
Several minimal models describing such self-organizing phenomenon
have been derived \cite{VicCziBenCohSho95, GreCha04,
CouKraFraLev05}. Most of these models include three basic effects:
short-range repulsion, long-range attraction, and reorientation or
alignment, in various ways, see \cite{HW} and particular
applications to birds \cite{HCH09} and fish \cite{BTTYB,BEBSVPSS}.

We first focus on populations of individuals driven by
self-propelling forces and pairwise attractive and repulsive
interaction \cite{LevRapCoh00, DorChuBerCha06}. We consider
self-propelled particles with Rayleigh friction
\cite{ChuHuaDorBer07, ChuDorMarBerCha07, CarDorPan09,UAB25},
leading to the Vlasov equation in $d=2,3$ dimensions:
\begin{equation}
\label{Equ1} \partial _t \fe + v \cdot \nabla _x \fe + a ^\eps
(t,x) \cdot \nabla _v \fe + \frac{1}{\eps} \Divv\{\fe \abv\}=
0,\;\;(t,x,v) \in \R_+ \times \R^d \times \R^d
\end{equation}
where $\fe = \fe (t,x,v) \geq 0$ represents the particle density
in the phase space $(x,v) \in \R^d \times \R^d$ at any time $t \in
\R_+$, $a ^\eps $ stands for the acceleration
\[
a^\eps (t,\cdot) = - \nabla _x U \star \rho ^\eps  (t, \cdot
),\;\;\rho ^\eps (t, \cdot ) = \intv{\fe (t, \cdot,
v)\;\mathrm{d}v}\, ,
\]
and $U$ is the pairwise interaction potential modelling the
repelling and attractive effects. Here, the propulsion and
friction forces coefficients $\alpha ^\eps =
\frac{\alpha}{\eps}>0$, $\beta ^\eps = \frac{\beta}{\eps}
>0$ are scaled in such a way that for $\eps\to 0$ particles will tend
to move with asymptotic speed $\sqrt{\tfrac{\alpha}\beta}$. These
models have been shown to produce complicated dynamics and
patterns such as mills, double mills, flocks and clumps, see
\cite{DorChuBerCha06}. Assuming that all individuals move with
constant speed also leads to spatial aggregation, patterns, and
collective motion \cite{CziStaVic97, EbeErd03}.

Another source of models arises from introducing alignment at the
modelling stage. A popular choice in the last years to include
this effect is the Cucker-Smale reorientation procedure
\cite{CS2}. Each individual in the group adjust their relative
velocity by averaging with all the others. This velocity averaging
is weighted in such a way that closer individuals in space have
more influence than further ones. The continuum kinetic version of
them leads to Vlasov-like models of the form \eqref{Equ1} in which
the acceleration is of the form
\[
a^\eps (t,\cdot) = - H \star f^\eps  (t, \cdot )\, ,
\]
where $\star$ stands for the $(x,v)$-convolution, abusing a bit on
the notation, with the nonnegative interaction kernel
$H:\R^{2d}\longrightarrow \R^d$. In the original Cucker-Smale
work, the interaction is modelled by $H(x,v)=h(x)v$, with the
weight function $h$ being a decreasing radial nonnegative
function. We refer to the extensive literature in this model for
further details \cite{HT08,HL08,CFRT10,review,MT11}.

In this work, we will consider the Vlasov equation \eqref{Equ1}
where the acceleration includes the three basic effects discussed
above, and then takes the form:
\begin{equation}\label{accel}
a^\eps (t,\cdot) = - \nabla _x U \star \rho ^\eps  (t, \cdot ) - H
\star f^\eps  (t, \cdot )\, .
\end{equation}
We will assume that the interaction potential $U\in C^2_b(\R^d)$,
$U$ bounded continuous with bounded continuous derivatives up to
second order, and $H(x,v)=h(x)v$ with $h\in C^1_b(\R^d)$ and
nonnegative. Under these assumptions the model
\eqref{Equ1}-\eqref{accel} can be rigorously derived as a
mean-field limit \cite{Neu77, BraHep77, Dob79,CCR10,BCC11} from
the particle systems introduced in \cite{DorChuBerCha06,CS2}.

We will first study in detail the linear problem, assuming that
the acceleration $a = a(t,x)$ is a given global-in-time bounded
smooth field. We investigate the regime $\eps \searrow 0$, that is
the case when the propulsion and friction forces dominate the
potential interaction between particles. At least formally we have
\begin{equation}
\label{EquAnsatz} \fe = f + \eps \fo + \eps ^2 f ^{(2)} + ...
\end{equation}
where
\begin{equation}
\label{Equ2} \Divv\{f \abv    \} = 0
\end{equation}
\begin{equation}
\label{Equ3} \partial _t f + \Divx (fv) + \Divv (f a(t,x)) +
\Divv\{\fo \abv    \} = 0\,,
\end{equation}
up to first order. Therefore, to characterize the zeroth order
term in the expansion we need naturally to work with solutions
whose support lies on the sphere of radius $r :=
\sqrt{\alpha/\beta}$ denoted by $r\sphere$ with $\sphere = \{v\in
\R^d : |v| = 1\}$. In turn, we need to work with measure solutions
to \eqref{Equ2} which makes natural to set as functional space the
set of nonnegative bounded Radon measures on $\R^d\times\R^d$
denoted by ${\cal M}_b ^+ (\R^d\times\R^d)$. We will be looking at
solutions to \eqref{Equ1} which are typically continuous curves in
the space ${\cal M}_b ^+ (\R^d\times\R^d)$ with a suitable notion
of continuity to be discussed later on. We will denote by
$\fe(t,x,v)\, \mathrm{d}(x,v)$ the integration against the measure
solution $\fe(t,x,v)$ of \eqref{Equ1} at time $t$. For the sake of
clarity, this is done independently of being the measure $\fe(t)$
absolutely continuous with respect to Lebesgue or not, i.e.,
having a $L^1(\R^d\times\R^d)$ density or not.

\begin{pro}\label{Kernel}
Assume that $(1+|v|^2)F \in {\cal M}_b ^+ (\R^d)$. Then $F$ is a
solution to \eqref{Equ2} if and only if $\supp F \subset \{0\}
\cup r \sphere$.
\end{pro}

The condition \eqref{Equ2} appears as a constraint, satisfied at
any time $t \in \R_+$. The time evolution of the dominant term $f$
in the Ansatz \eqref{EquAnsatz} will come by eliminating the
multiplier $\fo$ in \eqref{Equ3}, provided that $f$ verifies the
constraint \eqref{Equ2}. In other words we are allowed to use
those test functions $\psi (x,v)$ which remove the contribution of
the term $\Divv\{ \fo \abv \}$ {\it i.e.,}
\[
\intxv{\abv \cdot \nabla _v \psi \;\fo(t,x,v)\, \mathrm{d}(x,v) }
= 0.
\]
Therefore we need to investigate the invariants of the field $\abv
\cdot \nabla _v$. The admissible test functions are mainly those
depending on $x$ and $v/|v|, v \neq 0$. The characteristic flow
$(s,v) \to {\cal V}(s;v)$ associated to $\tfrac1\eps \abv \cdot
\nabla _v$
\[
\frac{\mathrm{d}{\cal V}}{\mathrm{d}s} = \frac1\eps
\abvs,\;\;{\cal V}(0;v) = v
\]
will play a crucial role in our study. It will be analyzed in
detail in Section \ref{LimMod}. Notice that the elements of $\A$
are the equilibria of $\abv \cdot \nabla _v $. It is easily seen
that the jacobian of this field
\[
\partial _v \{ \abv \} = (\alpha - \beta |v|^2 ) I - 2 \beta v \otimes v
\]
is negative on $r\sphere$, saying that $r\sphere$ are stable
equilibria. The point $0$ is unstable, $\partial _v \{ \abv \}
|_{v = 0}=\alpha I$. When $\eps \searrow 0$ the solutions
$(\fe)_\eps$ concentrate on $\xA$, leading to a limit curve of
measures even if $(\fe)_\eps$ were smooth solutions. We can
characterize the limit curve as solution of certain PDE whenever
our initial measure does not charge the unstable point $0$.

\begin{thm} \label{MainResult}
Assume that $a \in \litwoix{}$, $(1 + |v|^2) \fin \in \mbxv{}$,
$\supp \fin \subset \{(x,v) :|v|\geq r_0>0\}$. Then $(\fe)_\eps$
converges weakly $\star$ in $\litmbxv{}$ towards the solution of
the problem
\begin{equation}
\label{Equ22} \partial _t f + \Divx(fv) + \Divv \left \{f \imvv a \right \} = 0
\end{equation}
\begin{equation}
\label{Equ23} \Divv \{f \abv \} = 0
\end{equation}
with initial data $f(0) = \ave{\fin}$ defined by
$$
\intxv{\psi (x,v) \ave{\fin}(x,v)\, \mathrm{d}(x,v)} = \intxv{\psi
\left (x, r \vsv\right ) \fin(x,v)\, \mathrm{d}(x,v)}\,,
$$
for all $\psi \in \czcxv$.
\end{thm}

In the rest, we will refer to $\ave{\fin}$ as the projected
measure on the sphere of radius $r$ corresponding to $\fin$. Let
us point out that the previous result can be equivalently written
in spherical coordinates by saying that $f(t,x,\omega)$ is the
measure solution to the evolution equation on $(x,\omega)\in\R ^d
\times r \sphere$ given by
\begin{equation*}
\partial _t f + \Divx(f\omega) + \Divo \left \{f
\imoo a \right \} = 0 \,.
\end{equation*}
These results for the linear problem, when $a(t,x,v)$ is given, can
be generalized to the nonlinear counterparts where $a(t,x)$ is
given by \eqref{accel}. The main result of this work is (see Section \ref{MeaSol} for the definition of $\Po$):

\begin{thm} \label{MainResult2}
Assume that $U\in C^2_b(\R^d)$, $H(x,v)=h(x)v$ with $h\in
C^1_b(\R^d)$ nonnegative, $\fin \in \poxv{}$, $\supp \fin \subset
\{(x,v) :|x| \leq L_0, r_0\leq |v| \leq R_0\}$ with $0<r_0<r<R_0<\infty$. Then
for all $\delta>0$, the sequence $(\fe)_\eps$ converges in
$C([\delta,\infty);\poxv)$ towards the measure solution
$f(t,x,\omega)$ on $(x,\omega)\in\R ^d \times r \sphere$ of the
problem
\begin{equation}
\label{Equ22n} \partial _t f + \Divx(f\omega) - \Divo \left \{f
\imoo \left(\nabla_x U\star \rho + H\star f \right) \right \} = 0
\end{equation}
with initial data $f(0) = \ave{\fin}$. Moreover, if the initial
data $\fin$ is already compactly supported on $B_{L_0} \times r \sphere$, then
the convergence holds in $\cztpoxv$.
\end{thm}

Let us mention that the evolution problem \eqref{Equ22n} on $\R ^d
\times r \sphere$ was also proposed in the literature as the
continuum version \cite{DM08} of the Vicsek model
\cite{VicCziBenCohSho95,CouKraFraLev02} without diffusion for the
particular choice $U=0$ and $H(x,v)=h(x) v$ with $h(x)$ some local
averaging kernel. The original model in
\cite{VicCziBenCohSho95,CouKraFraLev02} also includes noise at the
particle level and was derived as the mean filed limit of some
stochastic particle systems in \cite{BCC12}. In fact, previous
particle systems have also been studied with noise in \cite{BCC11}
for the mean-field limit, in \cite{HLL09} for studying some
properties of the Cucker-Smale model with noise, and in
\cite{DFL10,FL11} for analyzing the phase transition in the Vicsek
model.

In the case of noise, getting accurate control on the particle
paths of the solutions is a complicated issue and thus, we are not
able to show the corresponding rigorous results to Theorems
\ref{MainResult} and \ref{MainResult2}. Nevertheless, we will
present a simplified formalism, which allows us to handle more
complicated problems to formally get the expected limit equations.
This approach was borrowed from the framework of the magnetic
confinement, where leading order charged particle densities have
to be computed after smoothing out the fluctuations which
correspond to the fast motion of particles around the magnetic
lines \cite{BosAsyAna, BosTraEquSin, BosGuiCen3D, BosNeg09}. We
apply this method to the following (linear or nonlinear) problem
\begin{equation}
\label{Equ31} \partial _t \fe + \Divx\{\fe v\} + \Divv \{ \fe a\} + \frac{1}{\eps} \Divv \{ \fe \abv \} = \Delta _v \fe
\end{equation}
with initial data $\fe (0) = \fin$ where the acceleration $a \in
\litwoix{}$ and $\fin \in \mbxv{}$. By applying the projection
operator $\ave{\cdot}$ to \eqref{Equ31}, we will show that the
limiting equation for the evolution of $f(t,x,\omega)$ on
$(x,\omega)\in\R ^d \times r \sphere$ is given by
\begin{equation} \label{Equ22Diff}
\partial _t f + \Divx(f\omega) + \Divo \left \{f \imoo a \right \} =
\Delta_\omega f
\end{equation}
where $\Delta_\omega$ is the Laplace-Beltrami operator on $r
\sphere$.

Our paper is organized as follows. In Section \ref{MeaSol} we
investigate the stability of the characteristic flows associated
to the perturbed fields $v \cdot \nabla _x + a \cdot \nabla _v +
\frac{1}{\eps} \abv \cdot \nabla _v $. The first limit result for
the linear problem (cf. Theorem \ref{MainResult}) is derived
rigorously in Section \ref{LimMod}. Section \ref{NLimMod} is
devoted to the proof of the main Theorem \ref{MainResult2}. The
new formalism to deal with the treatment of diffusion models is
presented in Section \ref{DiffMod}. The computations to show that
these models correspond to the Vicsek models, written in spherical
coordinates, are presented in the Appendix \ref{A}.


\section{Measure solutions}
\label{MeaSol}

\subsection{Preliminaries on mass transportation metrics and notations}
\label{prelim}

We recall some notations and result about mass transportation
distances that we will use in the sequel. For more details the
reader can refer to \cite{Vi1,CT}.

We denote by $\Po(\R^d)$ the space of probability measures on
$\R^d$ with finite first moment. We introduce the so-called
\emph{Monge-Kantorovich-Rubinstein distance} in $\Po(\R^d)$
defined by
\begin{equation*}
W_1(f,g) = \sup \left \{ \left |\int_{\R^d} \varphi(u)
(f(u)-g(u))\, \mathrm{d} u \right |, \varphi \in \lip(\R^d),
\lip(\varphi)\leq 1 \right \}
\end{equation*}
where $\lip(\R^d)$ denotes the set of Lipschitz functions on
$\R^d$ and $\lip(\varphi)$ the Lipschitz constant of a function
$\varphi$. Denoting by $\Lambda$ the set of transference plans
between the measures $f$ and $g$, i.e., probability measures in
the product space $\R^d \times \R^d$ with first and second
marginals $f$ and $g$ respectively
\[
f(y) = \int_{\R^d} \pi (y,z)\,\mathrm{d}z,\;\;g(z) = \int_{\R^d}
\pi (y,z)\,\mathrm{d}y
\]
then we have
\begin{equation*}
W_1(f, g) = \inf_{\pi\in\Lambda} \left\{ \int_{\R^d \times \R^d}
\vert y - z \vert \, \pi(y, z)\,\mathrm{d}(y,z) \right\}
\end{equation*}
by Kantorovich duality. $\Po(\R^d)$ endowed with this distance is
a complete metric space. Its properties are summarized below,
see\cite{Vi1}. 
\begin{pro}
\label{w2properties} The following properties of the distance
$W_1$ hold:
\begin{enumerate}
\item[1)] {\bf Optimal transference plan:} The infimum in the
definition of the distance $W_1$ is achieved. Any joint
probability measure $\pi_o$ satisfying:
$$
W_1(f, g) = \int_{\R^d \times \R^d} \vert y - z \vert \,
\mathrm{d}\pi_o(y, z)
$$
is called an optimal transference plan and it is generically non
unique for the $W_1$-distance.

\item[2)] {\bf Convergence of measures:} Given $\{f_k\}_{k\ge 1}$
and $f$ in $\Po(\R^d)$, the following two assertions are
equivalent:
\begin{itemize}
\item[a)] $W_1(f_k, f)$ tends to $0$ as $k$ goes to infinity.

\item[b)] $f_k$ tends to $f$ weakly $\star$  as measures as $k$
goes to infinity and
$$
\sup_{k\ge 1} \int_{\vert v \vert > R} \vert v \vert \, f_k(v) \,
\mathrm{d}v \to 0 \, \mbox{ as } \, R \to +\infty.
$$
\end{itemize}
\end{enumerate}
\end{pro}

Let us point out that if the sequence of measures is supported on
a common compact set, then the convergence in $W_1$-sense is
equivalent to standard weak-$\star$ convergence for bounded Radon
measures.

Finally, let us remark that all the models considered in this
paper preserve the total mass. After normalization we can consider
only solutions with total mass $1$ and therefore use the
Monge-Kantorovich-Rubinstein distance in $\Po (\R ^d \times \R
^d)$. From now on we assume that the initial conditions has total
mass $1$.

\subsection{Estimates on Characteristics}
In this section we investigate the linear Vlasov problem
\begin{equation}
\label{Equ10} \partial _t \fe + \Divx\{\fe v\} + \Divv \{ \fe a\}
+ \frac{1}{\eps} \Divv \{ \fe \abv \} = 0,\;\;(t,x,v) \in \R_+
\times \R^d \times \R^d
\end{equation}
\begin{equation}
\label{Equ11}
\fe (0) = \fin
\end{equation}
where $a \in \litwoix{}$ and $\fin \in \mbxv{}$.

\begin{defi}\label{DefMeaSol}
Assume that $a \in \litwoix{}$ and $\fin \in \mbxv{}$. We say that
$\fe \in \litmbxv{}$ is a measure solution of
\eqref{Equ10}-\eqref{Equ11} if for any test function $\varphi \in
\coctxv{}$ we have
\begin{align*}
\inttxv{\{\partial _t  + v \cdot \nabla _x + a \cdot \nabla _v +
\frac{1}{\eps} \abv \cdot & \nabla _v \}\varphi \fe(t,x,v)\,
\mathrm{d}(x,v)} \\
&+ \intxv{\varphi (0,x,v) \fin(x,v) \, \mathrm{d}(x,v) } = 0.
\end{align*}
\end{defi}
We introduce the characteristics of the field $v\cdot \nabla _x + a \cdot \nabla _v + \frac{1}{\eps} \abv \cdot \nabla _v $
\begin{equation*}
\frac{\mathrm{d}\Xe}{\mathrm{d}s} = \Ve(s),\;\;\frac{\mathrm{d}\Ve}{\mathrm{d}s} = a(s, \Xe(s)) + \frac{1}{\eps} \abves
\end{equation*}
\begin{equation*}
\Xe (s=0) = x,\;\;\Ve (s = 0) = v.
\end{equation*}
We will prove that $(\Xe, \Ve)$ are well defined for any $(s,x,v)
\in \R_+ \times \R^d \times \R^d$. Indeed, on any interval $[0,T]$
on which $(\Xe, \Ve)$ is well defined we get a bound
\[
\sup _{s \in [0,T]} \{|\Xe (s) | + |\Ve (s) |   \} < +\infty
\]
implying that the characteristics are global in positive time. For
that we write
\begin{equation}\label{charnew}
\frac12\frac{\mathrm{d}|\Ve|^2}{\mathrm{d}s} = a(s, \Xe(s))\cdot
\Ve (s) + \frac{1}{\eps} ( \alpha - \beta |\Ve (s) |^2) |\Ve
(s)|^2.
\end{equation}
and then, we get the differential inequality
\[
\frac{\mathrm{d}|\Ve |^2}{\mathrm{d}s} \leq 2\|a\|_{\linf} |\Ve
(s)| + \frac{2}{\eps} ( \alpha - \beta |\Ve (s) |^2) |\Ve (s)|^2
\]
for all $s\in [0,T]$, so that
\[
\sup _{s \in [0,T]}  |\Ve (s) | < +\infty,\;\;\sup _{s\in [0,T]} |\Xe (s) | \leq |x| + T \sup _{s\in [0,T]} |\Ve (s) | < +\infty.
\]
Once constructed the characteristics, it is easily seen how to
obtain a measure solution for the Vlasov problem
\eqref{Equ10}-\eqref{Equ11}. It reduces to push forward the
initial measure along the characteristics, see \cite{CCR10} for
instance.

\begin{pro}
For any $t \in \R_+$ we denote by $\fe (t)$ the measure given by
\begin{equation}\label{EquDefMea}
\intxv{\psi (x,v) \fe(t,x,v)\,\dxv} = \intxv{\psi((\Xe,
\Ve)(t;0,x,v))\fin(x,v)\,\dxv}\,,
\end{equation}
for all $\psi \in \czcxv$. Then the application  $t \to \fe (t)$,
denoted $\fin \#(\Xe, \Ve)(t;0,\cdot,\cdot)$ is the unique measure
solution of \eqref{Equ10}, \eqref{Equ11}, belongs to $\cztmbxv$
and satisfies
$$
\intxv{\fe (t,x,v)\,\dxv} = \intxv{\fin(x,v)\,\dxv}, t \in \R_+.
$$
\end{pro}
\begin{proof}
The arguments are straightforward and are left to the reader. We
only justify that $\fe \in \cztmbxv$ meaning that for any $\psi
\in \czcxv{}$ the application $t \to \intxv{\,\,\psi(x,v) \fe
(t,x,v)\;\dxv}$ is continuous. Choose $\psi\in \czcxv{}$. Then,
for any $0 \leq t_1 < t_2$ we have
\begin{align*}
\intxv{\psi(x,v)  \fe (t_2,x,v) &\,\dxv } - \intxv{\psi(x,v)  \fe
(t_1,x,v)\,\dxv } \\ &= \intxv{\left[\psi ((\Xe, \Ve )(t_2;t_1, x,
v)) - \psi (x,v)\right]\fe (t_1,x,v)\,\dxv}.
\end{align*}
Taking into account that $(\Xe, \Ve)$ are locally bounded (in
time, position, velocity) it is easily seen that for any compact
set $K \subset \R ^d \times \R^d$ there is a constant $C(K)$ such
that
\[
|\Xe (t_2; t_1, x, v) - x| + |\Ve (t_2; t_1, x, v) - v| \leq |t_2 - t_1 | C(K),\;\;(x,v) \in K.
\]
Our conclusion follows easily using the uniform continuity of $\psi$  and that
$\|\fe (t_1) \|_{{\cal M}_b} = \|\fin \|_{{\cal M}_b}$. Notice
also that the equality \eqref{EquDefMea} holds true for any
bounded continuous function $\psi$.
\end{proof}
We intend to study the behavior of $(\fe)_\eps$ when $\eps $
becomes small. This will require a more detailed analysis of the
characteristic flows $(\Xe, \Ve)$. The behavior of these
characteristics depends on the roots of functions like $A +
\frac{1}{\eps} (\alpha - \beta \rho ^2 ) \rho$, with $\rho \in
\R_+$, $A \in \R$.
\begin{pro}\label{NegA}
Assume that $A < 0$ and $ 0 < \eps < 2\alpha r /(|A|
3 \sqrt{3})$. Then the equation $\lae (\rho) := \eps A + (\alpha -
\beta \rho ^2 ) \rho = 0$ has two zeros on $\R_+$, denoted $\reo
(A), \ret (A)$, satisfying
\[
0 < \reo < \frac{r}{\sqrt{3}} < \ret < r
\]
and
\[
\lime \frac{\reo}{\eps} = \frac{|A|}{\alpha},\;\;\;\;\;\;\lime \frac{r - \ret}{\eps} = \frac{|A|}{2\alpha}
\]
where $r = \sqrt{\alpha/\beta}$.
\end{pro}
\begin{proof}
It is easily seen that the function $\lae$ increases on
$[0,r/\sqrt{3}]$ and decreases on $[r/\sqrt{3}, +\infty[$ with
change of sign on $[0,r/\sqrt{3}]$ and $[r/\sqrt{3}, r]$. We can
prove that $(\reo)_\eps, (\ret)_\eps$ are monotone with respect to
$\eps >0$. Take $0 < \eps < \teps <  2\alpha r /(|A| 3 \sqrt{3})$
and observe that $\lae > \lambda ^{\teps}$. In particular we have
\[
\lambda ^{\teps} (\reo) < \lae (\reo) = 0 = \lambda ^{\teps} (\rho _1 ^{\teps})
\]
implying $\reo < \rho _1 ^{\teps}$, since $\lambda ^{\teps}$ is
strictly increasing on $[0, r/\sqrt{3}]$. Similarly we have
\[
\lambda ^{\teps} (\ret) < \lae (\ret) = 0 < \lambda ^{\teps} (\rho _2 ^{\teps})
\]
and thus $\ret > \rho _2 ^{\teps}$, since $\lambda ^{\teps}$ is
strictly decreasing on $[r/\sqrt{3}, r]$. Passing to the limit in
$\lae (\rho _k ^\eps) = 0, k \in \{1,2\}$ it follows easily that
\[
\lime \reo = 0,\;\;\lime \ret = r.
\]
Moreover we can write
\begin{equation}
\alpha = \frac{\mathrm{d}}{\mathrm{d}\rho }\{(\alpha - \beta \rho ^2 ) \rho \} |_{\rho = 0} = \lime \frac{[\alpha - \beta (\reo)^2]\reo}{\reo} = - \lime \frac{\eps A}{\reo} \nonumber
\end{equation}
and
\begin{equation}
-2 \alpha = \frac{\mathrm{d}}{\mathrm{d}\rho }\{(\alpha - \beta \rho ^2 ) \rho \} |_{\rho = r} = \lime \frac{[\alpha - \beta (\ret)^2]\ret}{\ret - r} = - \lime \frac{\eps A}{\ret - r} \nonumber
\end{equation}
saying that
\[
\lime \frac{\reo}{\eps} = \frac{|A|}{\alpha},\;\;\lime \frac{r - \ret}{\eps} = \frac{|A|}{2\alpha}.
\]
\end{proof}
The case $A>0$ can be treated is a similar way and we obtain
\begin{pro} \label{PosA}
Assume that $A > 0$ and $ \eps >0$. Then the equation $\lae (\rho)
:= \eps A + (\alpha - \beta \rho ^2 ) \rho = 0$ has one zero on
$\R_+$, denoted $\reth (A)$, satisfying
\[
\reth >r,\;\;\lime \frac{\reth - r}{\eps} = \frac{|A|}{2\alpha}.
\]
\end{pro}
Using the sign of the function $\rho \to \eps \|a\|_{\linf{}} +
(\alpha - \beta \rho ^2 ) \rho$ we obtain the following bound for
the kinetic energy.

\begin{pro}\label{KinBou}
Assume that $a \in \litwoix{}$, $(1 + |v|^2) \fin \in \mbxv{}$ and
let us denote by $\fe$ the unique measure solution of
\eqref{Equ10}, \eqref{Equ11}. Then we have
\[
\left \|\intxv{\,|v|^2 \fe(\cdot,x,v)\,\dxv}\right \|_{\linf
(\R_+)} \leq \intxv{[(\reth)^2 + |v|^2] \fin(x,v)\,\dxv}.
\]
\end{pro}
\begin{proof}
We know that
\[
\frac{\mathrm{d}}{\mathrm{d}t} |\Ve|^2 \leq 2\|a\|_{\linf{}}
|\Ve(t)|+ \frac{2}{\eps} (\alpha - \beta |\Ve (t) |^2 ) |\Ve
(t)|^2=\frac{2}{\eps}|\Ve (t)|\lae (|\Ve (\overline{t})| ),\;\;t
\in \R_+.
\]
By comparison with the solutions of the autonomous differential
equation associated to the righthand side, we easily deduce that
\[
|\Ve (t;0,x,v)| \leq \max \{ |v|, \reth(\|a\|_{\linf{}})\}\,,
\]
for any $T \in \R_+, (x,v) \in \R ^d \times \R ^d$. This yields
the following bound for the kinetic energy
\begin{align*}
\intxv{|v|^2\fe (T,x,v)\,\dxv} &= \intxv{|\Ve (T;0,x,v)|^2
\fin(x,v)\,\dxv} \\
&\leq \intxv{[(\reth)^2 + |v|^2]
\fin(x,v)\,\dxv}.
\end{align*}
\end{proof}

The object of the next result is to establish the stability of
$\Ve$ around $|v| = r$. We will show that the characteristics
starting at points with velocities inside an annulus of length
proportional to $\eps$ around the sphere $r\sphere$ get trapped
there for all positive times for small $\eps$.

\begin{pro}
\label{RStab} Assume that $\eps \|a\|_{\linf{}} < 2\alpha r
/(3\sqrt{3})$ and that $\ret (-\|a\|_{\linf{}}) \leq |v| \leq
\reth (\|a\|_{\linf{}})$. Then, for any $(t,x) \in \R_+ \times
\R^d$ we have
\[
\ret (-\|a\|_{\linf{}}) \leq |\Ve(t;0,x,v)| \leq \reth (\|a\|_{\linf{}}).
\]
\end{pro}
\begin{proof}
As in previous proof, we know that
\[
\frac{\mathrm{d}}{\mathrm{d}t} |\Ve|^2 \leq \frac{2}{\eps}|\Ve
(t)|\lae (|\Ve (\overline{t})| ),\;\;t \in \R_+\,.
\]
By comparison with the constant solution $\reth$ to the autonomous
differential equation associated to the righthand side, we get
that $\sup _{t \in \R_+} |\Ve (t;0,x,v)| \leq  \reth$. Assume now
that there is $T>0$ such that $|\Ve (T) | < \ret$ and we are done
if we find a contradiction. Since $|\Ve (0) |= |v| \geq \ret$, we
can assume that $\min _{t \in [0,T]} |\Ve (t) | > \reo>0$ by time
continuity. Take now $\ot \in [0,T]$ a minimum point of $t \to
|\Ve (t)|$ on $[0,T]$. Obviously $\ot >0$ since
\[
|\Ve (\ot) | \leq |\Ve (T)| < \ret \leq |v| = |\Ve (0)|.
\]
By estimating from below in \eqref{charnew} and using that $\ot$
is a minimum point of $t \to |\Ve (t)|>0$ on $[0,T]$, we obtain
\[
0 \geq \frac{\mathrm{d}}{\mathrm{d}t} |\Ve (\ot)| \geq - \|a\|_{\linf{}} + \frac{(\alpha - \beta |\Ve (\ot)|^2)|\Ve (\ot)| }{\eps}
=\frac{\lae ( |\Ve (\ot)| )}{\eps}.
\]
But the function $\lae$ has negative sign on $[0,\reo] \cup [\ret,
+\infty[$. Since we know that $\min _{t \in [0,T]} |\Ve (t)| >
\reo$, it remains that
\[
\min _{t \in [0,T]} |\Ve (t)| = |\Ve (\ot)| \geq \ret
\]
which contradicts the assumption $|\Ve (T)| < \ret$.
\end{proof}

Let us see now what happens when the initial velocity is outside
$[\ret (-\|a\|_{\linf{}}), \reth (\|a\|_{\linf{}})]$. In
particular we prove that if initially $v \neq 0$, then $\Ve (t), t
\in \R_+$ remains away from $0$. We actually show that the
characteristics starting away from zero speed but inside the
sphere $r\sphere$ will increase their speed with respect to its
initial value while those starting with a speed outside the sphere
$r\sphere$ will decrease their speed with respect to its initial
value, all for sufficiently small $\eps$.

\begin{pro}
\label{ZeroStab} Consider $\eps >0$ such that $\eps \|a\|_{\linf{}} < 2\alpha r /(3\sqrt{3})$.\\
1. Assume that $\reo (- \|a\|_{\linf{}}) < |v| < \ret (-
\|a\|_{\linf{}})$. Then for any $(t,x) \in \R_+ ^\star \times \R
^d$ we have
\[
\reo (- \|a\|_{\linf{}}) < |v| < |\Ve (t;0,x,v)|\leq\reth ( \|a\|_{\linf{}}).
\]
2. Assume that $\reth ( \|a\|_{\linf{}}) < |v|$. Then for any
$(t,x) \in \R_+ ^\star \times \R^d$ we have
\[
\ret (- \|a\|_{\linf{}}) \leq |\Ve (t;0,x,v) | < |v|.
\]
\end{pro}
\begin{proof}
1. Notice that if $|\Ve (T;0,x,v)| = \ret$ for some $T>0$, then we
deduce by Proposition \ref{RStab} that $\ret \leq |\Ve (t) | \leq
\reth$ for any $t >T$ and thus $|\Ve (t;0,x,v) | \geq \ret > |v|,
t \geq T$. It remains to establish our statement for intervals
$[0,T]$ such that $|\Ve (t) | < \ret$ for any $t \in [0,T]$. We
are done if we prove that $t \to |\Ve (t)|$ is strictly increasing
on $[0,T]$. For any $\tau \in ]0,T]$ let us denote by $\ot$ a
maximum point of $t \to |\Ve (t)|>0$ on $[0,\tau]$. If $\ot \in
[0,\tau[$ we have $\frac{\mathrm{d}}{\mathrm{d}t} |\Ve (\ot)| \leq
0$ and thus
\[
0 \geq \frac{\mathrm{d}}{\mathrm{d}t} |\Ve (\ot)|\geq - \|a\|_{\linf{}} + \frac{(\alpha - \beta |\Ve (\ot)|^2)|\Ve (\ot)| }{\eps}
=\frac{\lae ( |\Ve (\ot)| )}{\eps}.
\]
By construction $|\Ve (\ot)| < \ret$ and moreover,
\[
|\Ve (\ot)| = \max _{[0,\tau]} |\Ve | \geq |v| > \reo\,,
\]
and thus, $\lae ( |\Ve (t)| )>0$ for all $t\in [0,T]$.
Consequently, we infer that $t \to |\Ve (t)|$ is strictly
increasing on $[0,T]$ since
\[
\frac{\mathrm{d}}{\mathrm{d}t} |\Ve (t)|\geq - \|a\|_{\linf{}} +
\frac{(\alpha - \beta |\Ve (t)|^2)|\Ve (t)| }{\eps} =\frac{\lae (
|\Ve (t)| )}{\eps} >0\,.
\]
Therefore we have $\ot = \tau$ saying that $|\Ve (\tau)| \geq |v|$
for any $\tau \in [0,T]$.

2. As before, it is sufficient to work on intervals $[0,T]$ such
that $|\Ve (t) | > \reth (\|a\|_{\linf{}})$ for any $t \in [0,T]$.
We are done if we prove that $t \to |\Ve (t)|$ is strictly
decreasing on $[0,T]$. We have for any $t \in [0,T]$
\[
\frac{\mathrm{d}}{\mathrm{d}t} |\Ve (t)|\leq  \|a\|_{\linf{}} + \frac{(\alpha - \beta |\Ve (t)|^2)|\Ve (t)| }{\eps}
=\frac{\lae ( |\Ve (t)| )}{\eps} <0
\]
where for the last inequality we have used $|\Ve (t) | > \reth, t \in [0,T]$.
\end{proof}


\section{The limit model}
\label{LimMod} We investigate now the stability of the family
$(\fe)_\eps$ when $\eps$ becomes small. After extraction of a
sequence $(\eps_k)_k$ converging to $0$ we can assume that
$(\fek)_k$ converges weakly $\star$ in $L^\infty(\R_+;{\cal M}_b
(\R^d \times \R^d))$, meaning that
\[
\limk \inttxv{\varphi (t,x,v) \fek (t,x,v)\,\dxv} =
\inttxv{\varphi (t,x,v) f (t,x,v)\,\dxv}
\]
for any $\varphi \in \lotczcxv{}$. Using the weak formulation of
\eqref{Equ10}-\eqref{Equ11} with test functions $\eta (t) \varphi
(x,v)$, $\eta \in C^1 _c (\R_+)$, $\varphi \in C^1 _c (\R^d \times
\R^d)$ one gets
\begin{align*}
 \inttxv{\{\eta ^{\;\prime} (t) \varphi + \eta (t) v \cdot \nabla _x \varphi + \eta (t) a \cdot \nabla _v \varphi   \}\fek(t,x,v)\,\dxv&}\\
+ \frac{1}{\eps _k} \inttxv{\eta (t) \abv \cdot \nabla _v \varphi \fek(t,x,v)\,\dxv &}  \\
= -\intxv{\eta (0) &\varphi (x,v) \fin(x,v)\,\dxv }.
\end{align*}
Multiplying by $\eps _k$ and passing to the limit for $k \to +\infty$ yields
\[
\inttxv{\eta (t) \abv \cdot \nabla _v \varphi f (t,x,v)\,\dxv} = 0
\]
and therefore one gets for any $t \in \R_+$ and $\varphi \in \cocxv{}$
\[
\intxv{\abv \cdot \nabla _v \varphi f (t,x,v)\,\dxv} = 0.
\]
Under the hypothesis $(1 + |v|^2) \fin \in \mbxv{}$ we deduce by
Proposition \ref{KinBou} that $( 1 + |v|^2) f(t) \in \mbxv{}$ and
therefore, applying the $(x,v)$ version of Proposition
\ref{Kernel} (whose proof is detailed in the sequel), we obtain
\[
\supp f(t) \subset \R^d \times (\A),\;\;t \in \R_+.
\]
The proof of Proposition \ref{Kernel} is based on the resolution of the adjoint problem
\[
- \abv \cdot \nabla _v \varphi = \psi (v),\;\;v \in \R^d
\]
for any smooth righthand side $\psi$ with compact support in $^c(\A)$.

\begin{proof} (of Proposition \ref{Kernel})
It is easily seen that for any $F \in \mbxv{}$, $\supp F \subset \A$ and any $\varphi \in \cocv{}$ we have
\[
\intv{\abv \cdot \nabla _v \varphi (v) F(v)\,\dv} = 0
\]
saying that $\Divv \{F \abv \} = 0$. Assume now that $\Divv \{F
\abv \} = 0$ for some $F \in \mbxv{}$ and let us prove that $\supp
F \subset \A$. We introduce the flow ${\cal V} = {\cal V}(s;v)$
given by
\begin{equation}
\label{Equ4} \frac{\mathrm{d}{\cal V}}{\mathrm{d}s} = ( \alpha - \beta |{\cal V} (s;v) |^2 ) {\cal V } (s;v),\;\;{\cal V}(0;v) = v.
\end{equation}
A direct computation shows that $\vsv$ are left invariant
\[
\abv \cdot \nabla _v \left ( \vsv \right ) = (\alpha - \beta |v|^2 )  \imvv \vsv = 0
\]
and therefore
\[
{\cal V} (s;v) = |{\cal V}(s;v)| \vsv,\;\;v \neq 0.
\]
Multiplying \eqref{Equ4} by ${\cal V}(s;v) / |{\cal V}(s;v)|$ yields
\[
\frac{\mathrm{d}}{\mathrm{d}s}|{\cal V}| = ( \alpha - \beta |{\cal V} (s;v) |^2 ) |{\cal V } (s;v)|
\]
whose solution is given by
\[
|{\cal V}(s;v)| = |v| \frac{r e ^{\alpha s}}{\sqrt{|v|^2 ( e ^{2\alpha s} - 1) + r^2}}
\]
Finally one gets
\begin{equation*}
{\cal V}(s;v) =  \frac{r e ^{\alpha s}}{\sqrt{|v|^2 ( e ^{2\alpha s} - 1) + r^2}}\;v,\;\;s \in ]S(v),+\infty[
\end{equation*}
with $S(v) = - \infty$ if $0 \leq |v| \leq r$ and $S(v) =
\frac{1}{2\alpha} \ln \left ( 1 - \frac{r^2}{|v|^2} \right ) < 0$
if $|v| > r$. Notice that the characteristics ${\cal V} (\cdot;v)$
are well defined on $\R_+$ for any $v \in \R^d$ and we have
\[
\lim _ {s \to +\infty} {\cal V}(s;v) = r \vsv\;\mbox{ if } v \neq 0,\;\;\lim _ {s \to +\infty} {\cal V}(s;v) =0\;\mbox{ if } v = 0
\]
and
\[
\lim _{s \searrow S(v)} |{\cal V}(s)| = 0\mbox{ if }0 \leq |v| < r,\;\lim _{s \searrow S(v)} |{\cal V}(s)| =r\mbox{ if } |v| = r,\;\lim _{s \searrow S(v)} |{\cal V}(s)| =+\infty\;\mbox{ if } |v| >r.
\]
Let us consider a $C^1$ function $\psi = \psi (v)$ with compact support in $^c (\A)$. We intend to construct a bounded $C^1$ function $\varphi = \varphi (v)$ such that
\begin{equation*}
- \abv \cdot \nabla _v \varphi = \psi (v),\;\;v \in
\R^d.
\end{equation*}
Obviously, if such a function exists, we may assume that $\varphi (0) = 0$. Motivated by the equality
\[
- \frac{\mathrm{d}}{\mathrm{d}s} \{\varphi ({\cal V}(s;v)) \}= \psi ({\cal V}(s;v)),\;\;0 \leq |v| < r,\;\;- \infty < s \leq 0
\]
and since we know that $\lim _{s \to - \infty} {\cal V} (s;v) = 0$ for any $0 \leq |v| < r$, we define
\begin{equation}
\label{Equ7} \varphi (v) = - \int _{-\infty} ^ 0 \psi ( {\cal V}(\tau; v))\;\mathrm{d}\tau,\;\;0 \leq |v| < r.
\end{equation}
Let us check that the function $\varphi$ in \eqref{Equ7} is well
defined and is $C^1$ in $|v|<r$. The key point is that $\psi $ has
compact support in $^c (\A)$ and therefore there are $0 < r_1 <
r_2 < r < r_3 < r_4 < +\infty$ such that $ \supp \psi \subset \{ v
\in \R ^d \;:\; r_1 \leq |v| \leq r_2 \} \cup \{ v \in \R^d \;:\;
r_3 \leq |v| \leq r_4\}. $ It is easily seen that $\tau \to |{\cal
V} (\tau; v)|$ is strictly increasing for any $0 < |v| < r$.
Therefore, for any $|v| \leq r_1$ we have $ |{\cal V} (\tau; v) |
\leq |{\cal V} (0; v) | = |v| \leq r_1,\;\;\tau \leq 0 $, implying
that
\[
\varphi (v) = - \int _{-\infty} ^ 0 \psi ({\cal V}(\tau; v))\;\mathrm{d}\tau = 0,\;\;0 \leq |v| \leq r_1.
\]
For any $v$ with $r_1 < |v| < r_2$ there are $\tau _1 < 0 < \tau _2$ such that
$
|{\cal V}(\tau _1; v)| = r_1 < r_2 = |{\cal V}(\tau _2; v)|.
$
The time interval between $\tau _1$ and $\tau _2$ comes easily by writing
\[
\frac{\frac{\mathrm{d}}{\mathrm{d}\tau}|{\cal V}(\tau) |}{(\alpha - \beta |{\cal V}(\tau)|^2)|{\cal V}(\tau) |}= 1
\]
implying that
\[
|\tau _2 | + |\tau _1 | = \tau _2 - \tau _1 = \int _{r_1} ^ {r_2} \frac{\mathrm{d}\rho}{(\alpha - \beta \rho ^2 ) \rho }.
\]
From the equality
\[
\varphi (v) = - \int _{-\infty} ^{\tau _1} \psi ({\cal
V}(\tau;v))\;\mathrm{d}\tau - \int _{\tau _1} ^0 \psi ({\cal
V}(\tau;v))\;\mathrm{d}\tau = - \int _{\tau _1} ^0 \psi ({\cal
V}(\tau;v))\;\mathrm{d}\tau\,,
\]
we deduce that
\begin{equation}
\label{Equ8} |\varphi (v) | \leq |\tau _1 | \; \|\psi \|_{C^0} \leq
\int _{r_1} ^ {r_2} \frac{\mathrm{d}\rho}{(\alpha - \beta \rho ^2 ) \rho }\; \|\psi \|_{C^0}.
\end{equation}
Assume now that $r_2 \leq |v| < r$. There is $\tau _2 \geq 0$ such
that $v = {\cal V} ( \tau_2 ; r_2 \vsv)$ and therefore
\begin{align*}
\varphi (v) & =  - \int _{-\infty} ^ 0 \psi ({\cal V}(\tau;v))\;\mathrm{d}\tau = - \int _{-\infty} ^ 0 \psi ({\cal V}(\tau + \tau _2;r_2 \vsv))\;\mathrm{d}\tau \\
& = - \int _{-\infty} ^ {-\tau _2} \psi ({\cal V}(\tau + \tau _2
;r_2 \vsv))\;\mathrm{d}\tau = - \int _{-\infty} ^ {0} \psi ({\cal
V}(\tau  ;r_2 \vsv))\;\mathrm{d}\tau = \varphi \left ( r_2 \vsv
\right).
\end{align*}
In particular, the restriction of $\varphi$ on $r_2 \leq |v| < r$
satisfies the same bound as in \eqref{Equ8}
\[
|\varphi (v) | \leq
\int _{r_1} ^ {r_2} \frac{\mathrm{d}\rho}{(\alpha - \beta \rho ^2 ) \rho }\; \|\psi \|_{C^0},\;\;r_2 \leq |v| < r.
\]
It is easily seen that $\varphi $ is $C^1$ on $0 \leq |v| < r$.
For that it is sufficient to consider $r_1 \leq |v| \leq r_2$.
Notice that
\[
\frac{\partial {\cal V}}{\partial v} (\tau; v) = \frac{|{\cal V}(\tau;v)|}{|v|} \left ( I - \frac{{\cal V}(\tau;v) \otimes {\cal V}(\tau;v)}{r^2} ( 1 - e ^ {-2\alpha \tau }  ) \right)
\]
and therefore the gradient of $\varphi$ remains bounded on $r_1
\leq |v| \leq r_2$
\[
\nabla _v \varphi (v) = - \int _{\tau _1} ^ 0 \frac{^ t \partial {\cal V}}{\partial v }(\tau; v) \nabla \psi ({\cal V}(\tau;v))\;\mathrm{d}\tau
\]
since on the interval $\tau \in [\tau _1, 0]$ we have $|{\cal
V}(\tau;v)| \in [r_1, |v|] \subset [r_1, r_2]$. Taking now as
definition for $|v| = r$
\[
\varphi (v) = \varphi \left ( r_2 \vsv \right )\,,
\]
we obtain a bounded $C^1$ function on $|v| \leq r$ satisfying
\[
- \abv \cdot \nabla _v \varphi = \psi (v),\;\;|\varphi (v) | \leq \int _{r_1} ^ {r_2} \frac{\mathrm{d}\rho}{(\alpha - \beta \rho ^2 ) \rho }\; \|\psi \|_{C^0},\;|v|\leq r.
\]
We proceed similarly in order to extend the above function for
$|v| > r$. We have for any $s>0$
\[
- \varphi ({\cal V}(s;v)) + \varphi (v) = \int _0 ^s \psi ({\cal V}(\tau;v))\;\mathrm{d}\tau,\;\;|v|> r.
\]
As $\lim _{s \to +\infty} {\cal V}(s;v) = r \vsv$ we must take
$$
\varphi (v) = \lim _{s \to +\infty}\left \{\varphi ( {\cal
V}(s;v)) + \int _0 ^s \psi ({\cal V}(\tau;v))\;\mathrm{d}\tau
\right \} = \varphi \left (r\vsv   \right ) + \int _0 ^{+\infty}
\psi ({\cal V}(\tau;v))\;\mathrm{d}\tau,\;\;|v| >r.\nonumber
$$
Clearly, for any $|v| > r$ the function $\tau \to |{\cal V}(\tau;v)|$ is strictly decreasing. Therefore, for any $r < |v| \leq r_3$ we have
\[
\varphi (v) = \varphi \left (r\vsv   \right )= \varphi \left (r_2\vsv   \right )
\]
since $|{\cal V}(\tau;v)|\leq |v| \leq r_3$ and $\psi ({\cal V}(\tau;v)) = 0$, $\tau \geq 0$. If $r_3 < |v| < r_4$ let us consider $\tau _4 < 0 < \tau _3$ such that
$
|{\cal V}(\tau _3;v)| = r_3 < r_4 = |{\cal V}(\tau _4;v)|.
$
The time interval between $\tau _4$ and $\tau _3$ is given by
\[
|\tau _3 | + |\tau _4 | = \tau _3 - \tau _4 = \int _{r_4} ^ {r_3}
\frac{\mathrm{d}\rho}{(\alpha - \beta \rho ^2) \rho } < +\infty\,,
\]
and therefore one gets for $r_3 < |v| < r_4$
\begin{align}
 |\varphi (v) | &\leq \left | \varphi \left ( r \vsv \right ) \right | + \left |\int _0 ^{\tau _3} \!\!\!\!\psi ({\cal V}(\tau;v))\;\mathrm{d}\tau \right | \nonumber \\
 & \leq \left [ \int _{r_1} ^ {r_2} \frac{\mathrm{d}\rho}{(\alpha - \beta \rho ^2) \rho } + \int _{r_4} ^ {r_3} \frac{\mathrm{d}\rho}{(\alpha - \beta \rho ^2) \rho } \right ] \|\psi \|_{C^0}.\label{Equ9}
\end{align}
Consider now $|v|\geq r_4$. There is $\tau _4 \geq 0$ such that $r_4 \vsv = {\cal V} (\tau_4; v)$ implying that
\begin{align*}
\varphi (v) & = \varphi \left ( r \vsv  \right ) + \int _0 ^{+\infty} \psi ({\cal V} (\tau; v)) \;\mathrm{d}\tau =
\varphi \left ( r \vsv  \right ) + \int _{\tau _4} ^{+\infty} \psi ({\cal V} (\tau; v)) \;\mathrm{d}\tau \\
& = \varphi \left ( r \vsv  \right ) + \int _0 ^{+\infty} \psi ({\cal V} (\tau; {\cal V}(\tau _4;v))) \;\mathrm{d}\tau
= \varphi \left ( r \vsv  \right ) + \int _0 ^{+\infty} \psi ({\cal V} (\tau; r_4 \vsv)) \;\mathrm{d}\tau \\
& = \varphi \left ( r_4 \vsv  \right ).
\end{align*}
We deduce that the restriction of $\varphi $ on $\{v :|v| \geq
r_4\}$ satisfies the same bound as in \eqref{Equ9}. Moreover the
function $\varphi $ is $C^1$ on $\{v:|v|\geq r\}$, with bounded
derivatives. Indeed, it is sufficient to consider only the case
$r_3 \leq |v| \leq r_4$, observing that
\begin{eqnarray}
\nabla _v \varphi (v)  =  \frac{r_2}{|v|} \imvv \nabla _v \varphi \left ( r_2 \vsv  \right ) + \int _0 ^{\tau _3} \frac{^t \partial {\cal V}}{\partial v }(\tau;v)\nabla \psi ({\cal V}(\tau;v)) \;\mathrm{d}\tau \nonumber
\end{eqnarray}
\[
|{\cal V} (\tau; v)| \in [r_3, |v| ] \subset [r_3, r_4],\;\tau \in [0,\tau_3],\;\;|\tau _3| + |\tau _4| = \int _{r_4} ^ {r_3} \frac{\mathrm{d}\rho}{(\alpha - \beta \rho ^2) \rho } < +\infty.
\]
By construction we have $- \abv \cdot \nabla _v \varphi = \psi
(v)$, $|v| >r$.

Consider a $C^1$ decreasing function on $\R_+$ such that $\chi
|_{[0,1]} = 1, \chi _{[2,+\infty[} = 0$. We know that
\[
\intv{\abv \cdot \nabla _v \left \{ \varphi (v) \chi \left (
\frac{|v|}{R} \right ) \right \}\,F(v)\,\dv} = 0,\;\;R>0\,,
\]
saying that
\[
\intv{\chi \left ( \frac{|v|}{R} \right )\abv \cdot \nabla _v \varphi
 \;F(v)\,\dv} + \intv{(\alpha - \beta |v|^2) \varphi (v) \frac{|v|}{R} \chi ^{\;\prime} \left ( \frac{|v|}{R} \right ) \;F(v)\,\dv} = 0.
\]
Since $\varphi$ and $\psi = -  \abv \cdot \nabla _v \varphi $ are
bounded and $F$ has finite mass and kinetic energy, we can pass to
the limit for $R \to +\infty$, using the dominated convergence
theorem. We obtain for any $C^1$ function $\psi$, with compact
support in $^c(\A)$
\[
\intv{\psi (v) F(v)\,\dv} = - \intv{\abv \cdot \nabla _v \varphi\,
F(v)\,\dv} = 0.
\]
Actually the previous equality holds true for any continuous
function $\psi$ with compact support in $^c(\A)$, since
$\intv{F(v)\,\dv} < +\infty$, so that $\supp F \subset \A$.
\end{proof}

In order to obtain stability for $(\fek)_k$ we need to avoid the
unstable equilibrium $v = 0$. For that we assume that the initial
support is away from zero speed: there is $r_0
>0$ (eventually small, let us say $r_0 < r$) such that
\begin{equation}
\label{Equ20} \supp \fin \subset \{ (x,v)\in \R ^d \times \R^d
\;:\;|v| \geq r_0\}.
\end{equation}

\begin{pro}
\label{UnifSupp} Under the hypothesis \eqref{Equ20} we have for any $\eps >0$ small enough
\[
\supp \fe (t) \subset \{ (x,v)\in \R ^d \times \R^d \;:\;|v| \geq
r_0\},\;\;t \in \R_+.
\]
\end{pro}

\begin{proof}
Take $\eps >0$ such that $\eps \|a\|_{\linf{}} < 2\alpha r /(3
\sqrt{3})$ and $\reo (- \|a\|_{\linf{}}) < r_0$. For any
continuous function $\psi = \psi (x,v)$ with compact support in
$\R ^d \times \{v\;:\; |v| < r_0\}$ we have
\begin{align*}
\intxv{\psi(x,v) \fe (t,x,v)\,\dxv} & = \intxv{\psi (\Xe
(t;0,x,v),
\Ve (t;0,x,v))\fin(x,v)\,\dxv } \\
& = \intxv{\psi (\Xe (t;0,x,v), \Ve (t;0,x,v) ){\bf 1}_{\{|v| \geq
r_0 \}}\fin(x,v)\,\dxv}.
\end{align*}
But for any $|v| \geq r_0 > \reo$ we know by Proposition
\ref{ZeroStab} that $|\Ve (t;0,x,v)| > |v| \geq r_0$, implying
that $\psi (\Xe (t), \Ve (t)) = 0$. Therefore one gets $\int
_{\R^d \times \R^d}{\psi(x,v) \fe (t,x,v)\,\dxv} = 0$ saying that
$\supp \fe (t) \subset \{ (x,v):|v| \geq r_0\}$.
\end{proof}

We are ready now to establish the model satisfied by the limit
measure $f$. The idea is to use the weak formulation of
\eqref{Equ10}, \eqref{Equ11} with test functions which are
constant along the flow of $\abv \cdot \nabla _v$, in order to get
rid of the term in $\frac{1}{\eps}$. These functions are those
depending on $x$ and $\vsv$. Surely, the invariants $\vsv$ have no
continuous extensions in $v = 0$, but we will see that we can use
it, since our measures $\fe$ vanish around $v = 0$.

\begin{proof} (of Theorem \ref{MainResult})
We already know that $f$ satisfies \eqref{Equ23}. Actually, since
$\supp \fe (t) \subset \{(x,v):|v|\geq r_0\}, t \in \R_+, \eps
>0$, we deduce that $\supp f(t) \subset \{(x,v):|v| \geq r_0\}$
and finally $\supp f(t) \subset \R ^d \times r \sphere, t \in
\R_+$. We have to establish \eqref{Equ22} and find the initial
data. Consider a $C^1$ decreasing function $\chi $ on $\R_+$ such
that $\chi |_{[0,1]} = 1, \chi _{[2,+\infty[} = 0$. For any $\eta
= \eta (t) \in C^1_c (\R_+)$, $\varphi = \varphi (x,v) \in
\cocxv{}$ we construct the test function
\[
\theta (t,x,v) = \eta (t) \left [ 1 - \chi \left ( \frac{2|v|}{r_0}\right )  \right ] \varphi \left ( x, r\vsv \right ).
\]
Notice that $\theta $ is $C^1$ and $\theta = 0$ for $|v| \leq
\frac{r_0}{2}$. When applying the weak formulation of
\eqref{Equ10}-\eqref{Equ11} with $\theta$, the term in
$\frac{1}{\eps}$ vanishes. Indeed, we can write
\begin{align*}
\frac{1}{\eps}\inttxv{\eta (t) & \abv \cdot \nabla _v \left \{\left [ 1 - \chi \left ( \frac{2|v|}{r_0}\right )  \right ]\varphi \left ( x, r\vsv \right ) \right \}\fe(t,x,v)\,\dxv } \nonumber \\
& = \frac{1}{\eps} \int _{\R_+} \eta (t) \int _{|v|\geq r_0} \abv
\cdot \nabla _v \left \{ \varphi \left ( x, r\vsv \right ) \right
\}\fe(t,x,v)\,\dxv \;\mathrm{d}t = 0.\nonumber
\end{align*}
For the term containing $\partial _t \theta$ we obtain the following limit when $k \to +\infty$
\begin{align*}
T_1 ^k := \inttxv{\partial _t \theta \fek(t,x,v)\,\dxv} \to &\inttxv{\partial _t \theta f(t,x,v)\,\dxv} \\
& = \int _{\R_+} \eta ^{\;\prime} (t) \int _{|v|\geq r_0} \varphi \left ( x, r\vsv \right ) f(t,x,v)\,\dxv \;\mathrm{d}t \\
& = \int _{\R_+} \eta ^{\;\prime} (t) \int _{|v| = r} \varphi \left ( x, r\vsv \right ) f(t,x,v)\,\dxv \;\mathrm{d}t \\
& = \int _{\R_+} \eta ^{\;\prime} (t) \int _{|v| = r} \varphi \left ( x, v\right ) f(t,x,v)\,\dxv \;\mathrm{d}t \\
& = \inttxv{\partial _t ( \eta \varphi )f(t,x,v)\,\dxv}.
\end{align*}
Similarly, one gets
\begin{align*}
T_2 ^k := \inttxv{ v \cdot \nabla _x \theta \fek (t,x,v)\,\dxv} \to & \inttxv{v \cdot \nabla _x  \theta f (t,x,v)\,\dxv} \\
& = \inttxv{v \cdot \nabla _x ( \eta \varphi )
f(t,x,v)\,\dxv}.\nonumber
\end{align*}
For the term containing $a \cdot \nabla _v \theta$ notice that on the set $|v| \geq r_0$ we have
\[
a \cdot \nabla _v \theta = \eta (t) a \cdot \nabla _v \left \{ \varphi \left ( x, r\vsv \right )\right \} = \eta (t) \frac{r}{|v|}a \cdot \imvv (\nabla _v \varphi ) \left ( x, r\vsv \right )
\]
and therefore we obtain
\begin{align*}
T_3 ^k := \inttxv{& a \cdot \nabla _v \theta \fek (t,x,v)\,\dxv} \to \inttxv{a \cdot \nabla _v \theta f (t,x,v)\,\dxv} \nonumber \\
& = \int _{\R_+} \eta  (t) \int _{|v|\geq r_0} \frac{r}{|v|} \imvv a \cdot (\nabla _v \varphi ) \left ( x, r\vsv \right ) f(t,x,v)\,\dxv \;\mathrm{d}t \\
& = \inttxv{\;\;\;\imvv a \cdot \nabla _v (\eta \varphi )
f(t,x,v)\,\dxv}.\nonumber
\end{align*}
For treating the term involving the initial condition, we write
\begin{align*}
T_4 : = \intxv{\theta (0,x,v) \fin(x,v)\,\dxv } &= \intxv{\eta (0)
\varphi \left ( x, r \vsv \right ) \fin(x,v)\,\dxv } \\
&= \intxv{\eta (0) \varphi (x,v) \ave{\fin}(x,v)\,\dxv}.
\end{align*}
Passing to the limit for $k \to +\infty$ in the weak formulation
$T_1 ^ k + T_2 ^ k + T_3 ^ k + T_4 = 0$ yields the problem
\[
\partial _t f + \Divx\{f v \} + \Divv \left \{f \imvv a   \right \} = 0,\;\;f(0) = \ave{\fin}
\]
as desired.
\end{proof}

\begin{remark}
\label{ConstraintPropagation} The constraint \eqref{Equ23} is
propagated by the evolution equation \eqref{Equ22}. This comes by
the fact that the flow $(X,V)$ associated to the field $v \cdot
\nabla _x + \imvv a \cdot \nabla _v$ leaves invariant $\R^d \times
r\sphere$. Indeed, if $(X,V)$ solves
\[
\frac{\mathrm{d}X}{\mathrm{d}s} = V(s),\;\;\frac{\mathrm{d}V}{\mathrm{d}s} = \left (I - \frac{V(s) \otimes V(s)}{|V(s)|^2} \right ) a(s, X(s))
\]
\[
X(s;0,x,v) = x,\;\;V(s;0,x,v) = v \neq 0
\]
then
\[
\frac{1}{2}\frac{\mathrm{d}}{\mathrm{d}s}|V(s)|^2 = \left (I - \frac{V(s) \otimes V(s)}{|V(s)|^2} \right ) a(s, X(s)) \cdot V(s) = 0
\]
saying that $|V(s;0,x,v)| = |v|$ for any $(s,x) \in \R_+ \times
\R^d$. In particular, for any continuous function $\psi = \psi
(x,v)$ with compact support in $^c (\R ^d \times r\sphere)$ we
have
\begin{align*}
\intxv{\psi(x,v) f(s,x,v)\,\dxv} & = \intxv{\psi (X(s;0,x,v),
V(s;0,x,v))
\ave{\fin}(x,v)\,\dxv} \\
& = \int _{|v| = r} \psi (X(s;0,x,v), V(s;0,x,v))
\ave{\fin}(x,v)\,\dxv = 0
\end{align*}
since $\supp \ave{\fin} \subset \R ^d \times r\sphere$. Therefore
for any $s \in \R_+$ we have $\supp f(s) \subset \R ^d \times
r\sphere$ implying that $\Divv \{f(s)\abv \} = 0, s \in \R_+$.
\end{remark}
\begin{remark}
\label{Uni} By the uniqueness of the solution for \eqref{Equ22}
with initial data $\ave{\fin}$, we deduce that all the family
$(\fe)_\eps$ converges weakly $\star$ in $\litmbxv{}$.
\end{remark}


\section{The non linear problem}
\label{NLimMod}

Up to now we considered the stability of the linear problems
\eqref{Equ10}-\eqref{Equ11} for a given smooth field $a = a(t,x)
\in \litwoix{}$. We concentrate now on the non linear problem
\begin{equation}
\label{Equ41} \partial _t \fe + \Divx\{\fe v\} + \Divv \{ \fe
a^\eps\} + \frac{1}{\eps} \{ \fe \abv \}= 0,\;\;(t,x,v) \in \R_+
\times \R ^d \times \R ^d
\end{equation}
with $a^\eps = - \nabla _x U \star \rho ^\eps - H \star \fe.$ The
well posedness of the non linear equation \eqref{Equ41} comes by
fixed point arguments in suitable spaces of measures, and it has
been discussed in \cite{CCR10, BCC12} in the measure solution
framework. We summarize next the properties of the solutions
$(\fe)_{\eps >0}$.

\begin{pro}
\label{ExiUniNonLin} Assume $h \in C^1_b (\R^d), U \in C^2 _b (\R
^d)$ and $( 1 + |v|^2 ) \fin \in \mbxv{}$. For all $\eps >0$,
there is a unique solution $(\fe, a ^\eps) \in C(\R_+;\poxv{})
\times \litwoix{}$ to
\begin{equation}
\label{Equ43} \partial _t \fe + \Divx \{\fe v \} + \Divv \{\fe
\aeps\} + \frac{1}{\eps}\Divv \{ \fe \abv \} = 0,\;\;(t,x,v) \in
\R _+ \times \R ^d \times \R ^d
\end{equation}
\begin{equation}
\label{Equ44} \aeps = - \nabla _x U \star \int _{\R ^d} \fe \;\md
v - H \star \fe,\;\;H(x,v) = h(x)v
\end{equation}
with initial data $\fe (0) = \fin$, satisfying the uniform bounds
\[
\sup _{\eps >0, t \in \R_+} \intxv{|v|^2 \fe
(t,x,v)\;\dxv}<+\infty
\]
\[
\sup _{\eps >0} \|\aeps \|_{L^\infty(\R_+;L^\infty(\R^d))} = :A
<+\infty,\;\;\sup _{\eps >0} \|\nabla _x \aeps
\|_{L^\infty(\R_+;L^\infty(\R^d))} = :A_1 <+\infty.
\]
Moreover, if the initial condition satisfies
\begin{equation*}
\supp \fin \subset \{ (x,v) \in \R ^d \times \R ^d
\;:\;|x| \leq L_0, r_0 \leq |v| \leq R_0 \}
\end{equation*}
for some $L_0 >0, 0 < r_0 < r < R_0 < +\infty$, then for any $\eps >0$
small enough we have
\[
\supp \fe (t) \subset \{ (x,v) \in \R ^d \times \R ^d \;:\;|x|\leq L_0 + t R_0, r_0 \leq |v| \leq R_0 \},\;\;t \in \R_+.
\]
\end{pro}
\begin{proof}
Here, we only justify the uniform bounds in $\eps$, the rest is a
direct application of the results in \cite{CCR10, BCC12}. The
divergence form of \eqref{Equ43} guarantees the mass conservation
\[
\intxv{\,\fe (t,x,v)\;\dxv } = \intxv{\,\fin (x,v)\;\dxv},\;\;t
\in \R_+.
\]
Notice that the term $- \Divv \{ \fe H \star \fe\}$ balances the
momentum
$$
\int _{\R^{2d}}{\!\! v \Divv \{ \fe H \star \fe\}\;\dxv} = \int
_{\R^{4d}}{\!\! h(x-x^{\prime}) (v ^{\prime}-v) \fe (t, x^\prime,
v^\prime)\fe (t,x,v)\;\md (x^\prime, v^\prime)}\dxv = 0
$$
and decreases the kinetic energy
\begin{align*}
\int _{\R^{2d}}{\!\!\!|v|^2 \Divv \{ \fe H \star \fe\}\;\dxv} &= 2\int _{\R^{4d}}{{\!\!\!\!h(x-x^{\prime}) (v ^{\prime}-v) \cdot v \fe (t, x^\prime, v^\prime)\fe (t,x,v)\;\md (x^\prime, v^\prime)}\dxv}  \\
& =  - \int _{\R^{4d}}{{\!\!h(x-x^{\prime}) |v - v ^{\prime}|^2
\fe (t, x^\prime, v^\prime)\fe (t,x,v)\;\md (x^\prime,
v^\prime)}\dxv}.
\end{align*}
In particular, as $|v|^2 \fin \in \mbxv{}$, then the kinetic
energy $\int _{\R^{2d}}{|v|^2 \fe (t,x,v)\;\dxv}$ remains bounded,
uniformly in time $t \in \R_+$ and $\eps >0$. Indeed, using the
continuity equation one gets
\[
\intxv{v \cdot (\nabla _x U \star \rho ^\eps ) \fe (t,x,v)\;\dxv}
= \frac{1}{2}\frac{\md}{\md t} \int _{\R^d}{(U \star \rho ^\eps
(t))(x) \rho ^\eps (t,x)\;\md x}
\]
and after multiplying \eqref{Equ43} by $\frac{|v|^2}{2}$ together
with \eqref{Equ44}, we obtain
\begin{align}
\frac{\md }{\md t} \intxv{&\,\,\left (\frac{|v|^2}{2} +\frac{U
\star \rho ^\eps}{2}   \right ) \fe (t,x,v) \;\dxv }
- \frac{1}{\eps}\intxv{(\alpha |v|^2 - \beta |v|^4) \fe (t,x,v) \;\dxv} \nonumber \\
& = - \frac{1}{2}\int_{\R^{4d }}h(x-x^{\prime}) |v - v
^{\prime}|^2 \fe (t, x^\prime, v^\prime)\fe (t,x,v)\;\md
(x^\prime, v^\prime)\dxv\leq 0\,.\label{energy}
\end{align}
Consider now $t ^\eps $ a maximum point on $[0,T], T>0$, of the
total energy
\[
W^\eps (t) = \intxv{\left (\frac{|v|^2}{2} + \frac{U \star \rho
^\eps}{2}   \right ) \fe (t,x,v) \;\dxv },\;\;t \in [0,T].
\]
If $t^\eps = 0$ then it is easily seen that for any $t \in [0,T]$
\[
\intxv{\frac{|v|^2}{2} \fe (t,x,v) \;\dxv} \leq
\intxv{\frac{|v|^2}{2} \fin (x,v) \;\dxv} + \|U \| _{\linf} \left
( \intxv{\fin \;\dxv}\right ) ^2.
\]
If $t ^\eps \in ]0,T]$ then $\frac{\md }{\md t} W ^\eps (t^\eps)
\geq 0$ implying from \eqref{energy} by moment interpolation in
$v$ that
\[
\sup _{\eps >0, T>0} \intxv{(1 + |v|^4) \fe (t ^\eps ,x,v)
\;\dxv } <+\infty
\]
and thus the inequality $W^\eps (t) \leq W ^\eps (t ^\eps), t \in
[0,T]$ yields
\begin{align*}
\sup _{\eps >0, t \in [0,T]} \intxv{\frac{|v|^2}{2} \fe (t,x,v) \;\dxv} \leq &  \sup _{\eps >0, T>0} \intxv{\frac{|v|^2}{2} \fe (t ^\eps ,x,v) \;\dxv } \\
& + \|U \|_{\linf} \left ( \intxv{\fin \;\dxv}\right ) ^2
<+\infty.
\end{align*}
Therefore the kinetic energy remains bounded on $[0,T]$, uniformly
with respect to $\eps >0$, and the bound does not depend on $T>0$.
The uniform bounds for $\aeps$ come immediately by convolution
with $\nabla _x U$ and $H$, thanks to the uniform estimate
\[
\sup _{\eps >0, t \in \R_+} \intxv{|v| \fe (t,x,v)} < +\infty.
\]
We analyze the support of $(\fe )_{\eps >0}$. Take $\eps >0$ small
enough such that $\eps A < 2 \alpha r /(3\sqrt{3})$ and $\reo (-A)
< r_0,\; \reth (A) < R_0$. By Proposition \ref{UnifSupp} we already
know that
\[
\supp \fe (t) \subset \{ (x,v) \in \R ^d \times \R ^d \;:\;|v |
\geq r_0\},\;\;t \in \R_+.
\]
For any continuous function $\psi = \psi (x,v)$ with compact
support in $\R ^d \times \{ v\in \R ^d\;:\;|v| > R_0\}$ we have
\begin{align*}
\intxv{\psi(x,v) \fe(t,x,v) \;\dxv} & =
\intxv{\psi (\Xe (t), \Ve (t)) \fin(x,v) \;\dxv } \\
& = \intxv{\psi (\Xe (t), \Ve (t) ) \ind{\{r_0 \leq |v| \leq
R_0\}}\fin(x,v) \;\dxv}.
\end{align*}
We distinguish several cases:\\
1. If $r_0 \leq |v| < \ret (-A)$ we deduce by Proposition
\ref{ZeroStab} that
$
|v| < |\Ve (t;0,x,v)| \leq \reth (A) < R_0,\;\;t \in \R_+, \eps
>0.
$

2. If $\ret (-A) \leq |v| \leq \reth (A)$ we obtain by Proposition
\ref{RStab} that
$
\ret (-A) \leq |\Ve (t;0,x,v)| \leq \reth (A) < R_0,\;\;t \in
\R_+, \eps >0.
$

3. If $\reth (A) < |v| \leq R_0$ one gets thanks to Proposition
\ref{ZeroStab}
$
\ret (-A) \leq |\Ve (t;0,x,v) | < |v| \leq R_0.
$

In all cases $(\Xe, \Ve )(t;0,x,v)$ remains outside the support of
$\psi$, implying that
\[
\intxv{\psi (x,v) \fe (t,x,v)\;\dxv} = 0.
\]
Thus for any $t \in \R_+$ and $\eps >0$ small enough one gets
\[
\supp \fe (t) \subset \{ (x,v) \in \R ^d \times \R ^d \;:\;r_0
\leq |v | \leq R_0\}.
\]
Consider $\theta \in C^1 (\R)$ non decreasing, verifying $\theta
(u) = 0$ if $u \leq 0$, $\theta (u) >0$ if $u>0$. Applying the
weak formulation of \eqref{Equ43}-\eqref{Equ44} with the test
function $\theta (|x| - L_0 - t R_0)$ yields
\begin{align*}
\intxv{\theta (|x| - L_0 -& t R_0) \fe (t,x,v)\;\dxv}  = \intxv{\theta (|x| - L_0 )\fin(x,v)\;\dxv} \\
& + \int _0 ^t \intxv{\theta ^{\prime}(|x| - L_0 - s R_0) \left (
v \cdot \frac{x}{|x|} - R_0\right )\fe (s,x,v)\;\dxv} \md s \leq 0
\end{align*}
implying that  $\supp \fe (t) \subset \{ (x,v) \in \R ^d \times \R ^d :|x| \leq L_0 + t R_0\}, t \in \R_+$.
\end{proof}

\

The uniform bound for the total mass allows us to extract a
sequence $(\eps _k)_k \subset \R_+ ^\star$ convergent to $0$ such
that $(\fek)_k$ converges weakly $\star$ in $\litmbxv{}$. The
treatment of the non linear term requires a little bit more, that
is convergence in $C(\R_+;\poxv{})$ or at least in $C([\delta,
+\infty[;\poxv{})$ for any $\delta >0$. The key argument for
establishing that is emphasized by the lemma
\begin{lemma}
\label{TimeEstimate} Consider $\eps >0$ small enough.\\
1. For any $(x,v)\in \R ^d \times \R ^d$ with $r_0 \leq |v| < \ret
(-A) - \eps$, the first time $t^\eps _1 = t ^\eps _1 (x,v)$ such
that $|\Ve (t^\eps _1;0,x,v) | = \ret (-A) - \eps$ satisfies
\[
t ^\eps _1 \leq \frac{\eps}{2\beta r_0 ^2} \ln \left ( \frac{r -
r_0 }{\eps} \right ).
\]
2. For any $(x,v)\in \R ^d \times \R ^d$ with $\reth (A) + \eps  <
|v| \leq R_0$, the first time $t^\eps _2 = t ^\eps _2 (x,v)$ such
that $|\Ve (t^\eps _2;0,x,v) | = \reth (A) + \eps$ satisfies
\[
t ^\eps _2 \leq \frac{\eps}{2\beta r ^2} \ln \left ( \frac{R_0 -
r}{\eps} \right ).
\]
\end{lemma}
\begin{proof}
1. During the time $[0,t ^\eps _1]$ the velocity modulus $|\Ve
(t)|$ remains in $[r_0, \ret (-A) - \eps] \subset [\reo (-A), \ret
(-A)]$ and we can write for any $t \in [0, t ^\eps _1]$
\[
\frac{\eps \frac{\md |\Ve |}{\md t }}{- \eps A + (\alpha - \beta
|\Ve (t) |^2 ) \;|\Ve (t) |}\geq \frac{\frac{\md |\Ve |}{\md t
}}{\aeps (t, \Xe (t)) \cdot \frac{\Ve (t)}{|\Ve (t)|} +
\frac{1}{\eps}(\alpha - \beta |\Ve (t) |^2 ) \;|\Ve (t) |} = 1
\]
since $- \eps A + (\alpha - \beta u^2)u$ is positive for $u \in
[\reo (-A), \ret (-A)]$. Integrating with respect to $t \in [0, t
^\eps _1]$ yields
\[
t ^\eps _1 (x,v) \leq \eps \int _{|v|} ^{\ret (-A) - \eps }
\frac{\md u }{- \eps A + (\alpha - \beta u ^2 ) u } \leq \eps \int
_{r_0} ^{\ret (-A) - \eps } \frac{\md u }{- \eps A + (\alpha -
\beta u ^2 ) u }.
\]
Recall that $\ret (-A)$ is one of the roots of $u \to - \eps A +
(\alpha - \beta u ^2 ) u$ and therefore a direct computation lead
to
\[
- \eps A + (\alpha - \beta u ^2 ) u  = \beta (\ret - u ) [u ^2 + u
\ret + (\ret ) ^2 - r^2] \geq 2 \beta r_0 ^2 ( \ret - u), \;\;u
\in [r_0, \ret],\;\eps \;\mbox{small enough }
\]
implying that
\[
t ^\eps _1 (x,v) \leq \frac{\eps}{2\beta r_0 ^2} \int _{r_0}
^{\ret - \eps} \frac{\md u }{ \ret - u} = \frac{\eps}{2\beta r_0
^2}\ln \left (\frac{\ret - r_0}{\eps}   \right ) \leq
\frac{\eps}{2\beta r_0 ^2}\ln \left (\frac{r - r_0}{\eps}   \right
).
\]
2. During the time $[0,t ^\eps _2]$ the velocity modulus $|\Ve
(t)|$ remains in $[\reth (A) + \eps, R_0] \subset [\reth(A),
+\infty[$ and we can write for any $t \in [0, t ^\eps _2]$
\[
\frac{\eps \frac{\md |\Ve |}{\md t }}{ \eps A + (\alpha - \beta
|\Ve (t) |^2 ) \;|\Ve (t) |}\geq \frac{\frac{\md |\Ve |}{\md t
}}{\aeps (t, \Xe (t)) \cdot \frac{\Ve (t)}{|\Ve (t)|} +
\frac{1}{\eps}(\alpha - \beta |\Ve (t) |^2 ) \;|\Ve (t) |} = 1
\]
since $ \eps A + (\alpha - \beta u^2)u$ is negative for $u \in
[\reth (A), +\infty[$. Integrating with respect to $t \in [0, t
^\eps _2]$ yields
\[
t ^\eps _2 (x,v) \leq \eps \int _{|v|} ^{\reth (A) + \eps }
\frac{\md u }{ \eps A + (\alpha - \beta u ^2 ) u } \leq \eps \int
_{R_0} ^{\reth (A) + \eps } \frac{\md u }{ \eps A + (\alpha - \beta
u ^2 ) u }.
\]
By direct computation we obtain
\[
 \eps A + (\alpha - \beta u ^2 ) u  = - \beta (u - \reth ) [u ^2 + u \reth + (\reth ) ^2 - r^2] \leq -2 \beta r ^2 ( u - \reth ), \;\;u \geq \reth,\;\eps \;\mbox{small enough }
\]
implying that
\[
t ^\eps _2 (x,v) \leq \frac{\eps}{2\beta r ^2} \int _{\reth +
\eps} ^{R_0} \frac{\md u }{ u - \reth} = \frac{\eps}{2\beta r
^2}\ln \left (\frac{R _0 - \reth }{\eps}   \right ) \leq
\frac{\eps}{2\beta r ^2}\ln \left (\frac{R_0 - r}{\eps}   \right
).
\]
\end{proof}
We intend to apply Arzela-Ascoli theorem in $C(\R_+;\Po(\R ^d
\times \R ^d))$ in order to extract a convergent sequence
$(\fek)_k$ with $\limk \eps _k = 0$. We need to establish the
uniform equicontinuity of the family $(\fe)_{\eps >0}$. The
argument below is essentially similar to arguments in
\cite{CCR10}.

\begin{pro}
\label{UnifEquiCont} 1. If the initial data is well prepared {\it
i.e.,} $\supp \fin \subset \{ (x,v) \in \R ^d \times \R ^d\;:\;|x|\leq L_0,
|v| = r\}$ then there is a constant $C$ (not depending on $t \in
\R_+, \eps >0$) such that
\[
W_1 (\fe (t), \fe (s)) \leq C | t - s|,\;\;t, s \in \R_+, \eps >0.
\]
2. If $\supp \fin \subset \{ (x,v) \in \R ^d \times \R ^d\;:\;|x|\leq L_0, r_0
\leq |v| \leq R_0\}$ then there is a constant $C$ (not depending
on $t \in \R_+, \eps >0$) such that for any $\delta >0$ we can
find $\eps _\delta$ satisfying
\[
W_1 (\fe (t), \fe (s)) \leq C |t - s|,\;\; t,s \geq \delta,\;\;0 <
\eps < \eps _\delta.
\]
\end{pro}
\begin{proof}
1. Consider $\varphi = \varphi (x,v)$ a Lipschitz function on
$\R^d \times \R ^d$ with $\lip (\varphi ) \leq 1$. For any $t, s
\in \R_+, \eps >0$ we have
\begin{align*}
\left | \intxv{\!\!\!\varphi  ( \fe (t) - \fe (s))\dxv}\right |  &=  \left | \intxv{\!\!\!\{ \varphi (\Xe (t), \Ve (t)) - \varphi (\Xe (s), \Ve (s))\}\fin (x,v)\dxv}\right |  \\
& \leq  \intxv{\{ |\Xe (t) - \Xe (s)| + |\Ve (t) - \Ve (s)|\}
\ind{\{|v| = r\}} \fin \dxv}.
\end{align*}
Thanks to Proposition \ref{RStab} we have for any $(\tau, x, v)
\in \R_+ \times \R ^d \times r \sphere $
\[
\frac{\ret (-A) - r}{\eps} \leq \frac{|\Ve (\tau;0,x,v)| -
r}{\eps} \leq \frac{\reth (A) - r}{\eps}
\]
and it is easily seen, integrating the system of characteristics
between $s$ and $t$, that
\[
|\Xe (t;0,x,v) - \Xe (s;0,x,v) | = \left | \int _s ^ t \Ve
(\tau;0,x,v) \;\mathrm{d}\tau\right | \leq R_0 |t-s|
\]
and
\begin{align*}
\left | \Ve (t;0,x,v) - \Ve (s;0,x,v) \right | & \leq  \left | \int _s ^t \left \{  |a^\eps(\tau, \Xe  (\tau))| + \frac{|\alpha - \beta |\Ve (\tau) | ^2 | \;|\Ve (\tau) |}{\eps} \right \}\md \tau  \right |\nonumber \\
& \leq  |t -s | \left \{ A + \beta ( r + R_0) R_0 \max \left (
\frac{\reth (A) - r}{\eps}, \frac{r - \ret (-A) }{\eps} \right )
\right \}.
\end{align*}
Our conclusion comes immediately by Propositions \ref{NegA}, \ref{PosA}.\\
2. Consider $\delta >0$ and $\eps _\delta $ small enough such that
$\frac{\eps}{2\beta r_0 ^2} \ln \left ( \frac{r - r_0}{\eps}
\right ) < \delta$, $\frac{\eps}{2\beta r ^2} \ln \left (
\frac{R_0 - r}{\eps} \right ) < \delta$ for $0 < \eps < \eps
_\delta$. For any Lipschitz function $\varphi $ with $\lip
(\varphi ) \leq 1$ and any $t, s \geq \delta$ we have
\[
\left | \intxv{\!\!\!\!\varphi ( \fe (t) - \fe (s) ) \;\dxv}\right | \leq \!\!\intxv{\{|\Xe
(t) - \Xe (s) | + |\Ve (t) - \Ve (s)|  \} \ind{\{r_0 \leq |v| \leq
R_0\}} \fin \;\dxv}.
\]
For any $(\tau, x) \in \R_+ \times \R ^d$, $\ret (-A) - \eps \leq
|v| \leq \reth (A) + \eps$ we have by Propositions \ref{RStab},
\ref{ZeroStab}
\[
\ret (-A) - \eps \leq |\Ve (\tau;0,x,v) | \leq \reth (A) + \eps.
\]
The same conclusion holds true for any $\tau \geq \delta$, $x \in
\R^d$ and $|v| \in [r_0, \ret (-A) - \eps[ \cup ]\reth (A) + \eps,
R_0]$, thanks to Lemma \ref{TimeEstimate}, since $\delta > \max \{
t^\eps _1 (x,v), t^\eps _2 (x,v)\}$ (after a time $\delta$, the
velocity modulus $|\Ve (\tau;0,x,v)|$ is already in the set
$\{w\;:\;\ret (-A) - \eps < |w| < \reth (A) + \eps \}$). Our
statement follows as before, integrating the system of
characteristics between $s$ and $t$.
\end{proof}
Applying Arzela-Ascoli theorem, we deduce that there is a sequence
$(\eps _k)_k \subset \R _+ ^\star$, convergent to $0$ such that
\[
\limk W_1 (\fek (t), f(t)) = 0 \mbox{ uniformly for } t \in
[0,T],\;\; T>0
\]
for some $f \in C(\R_+;\Po (\R ^d \times \R ^d))$ if $\supp \fin
\subset \{(x,v)\in \R ^d \times \R ^d\;:\;|x| \leq L_0,|v| = r\}$ and
\[
\limk W_1 (\fek (t), f(t)) = 0 \mbox{ uniformly for } t \in
[\delta ,T],\;\; T>\delta >0
\]
for some $f \in C(\R_+ ^\star;\Po (\R ^d \times \R ^d))$ if $\supp
\fin \subset \{(x,v)\in \R ^d \times \R ^d\;:\;|x|\leq L_0, r_0 \leq |v| \leq
R_0\}$. It is easily seen that if the initial condition is well
prepared then there is a constant $C$ cf. Proposition
\ref{UnifEquiCont} such that
$
W_1 (f(t), f(s)) \leq C |t -s |,\;\;t, s \in \R_+.
$
The same is true for not prepared initial
conditions $\fin$. Take $\delta >0$ and $\eps _\delta$ as in
Proposition \ref{UnifEquiCont}. For any $0 < \eps < \eps _\delta$
we have
$
W_1 (\fe(t), \fe(s)) \leq C |t -s |,\;\;t, s \geq \delta.
$
For $k$ large enough we have $\eps _k < \eps _\delta$ and
therefore
$
W_1 (\fek(t), \fek(s)) \leq C |t -s |,\;\;t, s \geq \delta.
$
Passing to the limit as $k$ goes to infinity yields
$
W_1 (f(t), f(s)) \leq C |t -s |,\;\;t, s \geq \delta.
$
Since the constant $C$ does not depend on $\delta$ one gets
\[
W_1 (f(t), f(s)) \leq C |t -s |,\;\;t, s >0.
\]
In particular we deduce that $f$ has a limit as $t$ goes to $0$
since $( \Po (\R ^d \times \R ^d), W_1)$ is a complete metric
space and therefore we can extend $f$ by continuity at $t = 0$.
The extended function, still denoted by $f$, belongs to
$C(\R_+;\Po (\R ^d \times \R ^d))$ and satisfies
\[
W_1 (f(t), f(s)) \leq C |t -s |,\;\;t, s \in \R_+.
\]
The above convergence allows us to handle the non linear terms. We
use the following standard argument \cite{Dob79,CCR10}.

\begin{lemma}
\label{NonLinTerm} Consider $f,g \in \Po (\R ^d \times \R ^d)$
compactly supported $ \supp f \cup \supp g \subset \{ (x,v) \in \R
^d \times \R ^d \;:\;|x|\leq L, |v| \leq R\}$, and let us consider
\[
a_f = - \nabla _x U \star \int _{\R ^d} f \;\md v - H \star
f,\;\;a_g = - \nabla _x U \star \int _{\R ^d} g \;\md v - H \star
g.
\]
Then we have
\[
\|a_f - a_g \|_{L^\infty ( \R ^3 \times B_R)} \leq \left
\{\|\nabla _x ^2 U \|_{\linf} + \left (\|h \|^2 _{\linf} + 4 R ^2
\|\nabla _x h \| ^2 _{\linf} \right ) ^{1/2}   \right \}W_1 (f,g)
\]
where $B_R$ stands for the closed ball in $\R ^d$ of center $0$
and radius $R$.
\end{lemma}
\begin{proof}
Take $\pi $ to be a optimal transportation plan between $f$ and
$g$. Then for any $x \in \R^d$ we have, using the marginals of
$\pi$
\begin{align*}
| (\nabla _x U \star f) (x) - (\nabla _x U \star g) (x) | &=  \left | \intxv{\nabla _x U ( x - x^\prime) \{ f(\xp, \vp) - g(\xp, \vp)\}\;\dxpvp} \right |  \\
& =  \left | \intxv{\intxv{ [\nabla _x U (x - \xp) - \nabla _x U (x - \xs)]\md \pi (\xp, \vp, \xs, \vs)}}  \right |  \\
& \leq  \|\nabla _x ^2 U \|_{\linf{}} \intxv{\intxv{|\xp - \xs | \;\md \pi (\xp, \vp, \xs, \vs)}} \\
& \leq  \|\nabla _x ^2 U \|_{\linf{}} W_1 (f,g).
\end{align*}
The estimate for $H \star f - H \star g$ follows similarly
observing that on the support of $\pi$, which is included in
$\{(\xp, \vp, \xs, \vs)\in \R ^{4d}\;:\; |\vp|\leq R, |\vs | \leq
R\}$ we have
\begin{align*}
|h(x- \xp) (v- \vp) - h(x- \xs) & (v- \vs)| \\ & \leq  |h(x- \xp)( \vs - \vp)| + |h(x- \xp)- h(x- \xs)| \;|v - \vs |  \\
& \leq  \left ( \|h \|^2 _{\linf} + 4 R ^2 \|\nabla _x h \| ^2
_{\linf}  \right ) ^{1/2}  \left ( |\xp - \xs | ^2 + |\vp - \vs
|^2  \right ) ^{1/2}.
\end{align*}
\end{proof}

We are ready now to prove Theorem \ref{MainResult2}.

\begin{proof}
(of Theorem \ref{MainResult2}) The arguments are the same as those
in the proof of Theorem \ref{MainResult} except for the treatment
of the non linear terms. We only concentrate on it. Consider
$(\fek )_k$ with $\limk \eps _k = 0$ such that $\limk W_1 (\fek
(t), f(t)) = 0$ uniformly for $t \in [0,T], T>0$ if $\supp \fin
\subset \{ (x,v) \;:\;|x|\leq L_0, |v| = r\}$ and $\limk W_1 (\fek
(t), f(t)) = 0$ uniformly for $t \in [\delta,T], T>\delta >0$ if
$\supp \fin \subset \{ (x,v) \;:\;|x|\leq L_0, r_0 \leq |v| \leq
R_0\}$ for some function $f \in C(\R_+;\Po (\R ^d \times \R ^d))$.
Thanks to Proposition \ref{Kernel} we deduce (for both prepared or
not initial data) that
\[
\supp f(t)  \subset \{ (x,v)\in \R ^d \times \R ^d \;:\;|v| =
r\},\;\;t>0.
\]
The previous statement holds also true at $t = 0$, by the
continuity of $f$. The time evolution for the limit $f$ comes by
using the particular test functions
\[
\theta (t,x,v) = \eta (t) \left [ 1 - \chi \left (
\frac{2|v|}{r_0}\right )  \right ] \varphi \left ( x, r\vsv \right
)
\]
with $\eta \in C^1_c (\R_+)$, $\varphi \in \cocxv{}$. From now on
we consider only the not prepared initial data case (the other
case is simpler). We recall the notation $\aeps = - \nabla _x U
\star \int _{\R^d} \fe \;\md v - H \star \fe $ and we introduce $a = - \nabla _x
U \star \int _{\R^d} f \;\md v - H \star f$. Since $f$ satisfies the same
bounds as $(\fe)_\eps$, we deduce that $\|a\|_{\linf{}} \leq A,
\|\nabla _x a\|_{\linf{}} \leq A_1$. For any $\delta >0$ we can
write
\begin{align}
\label{EquBil}
&\left |\inttxv{\left \{ \aek \cdot \nabla _v \theta \;\fek   -   a \cdot \nabla _v \theta \;f\right \}\dxv\!\!}  \right |
\leq  \left |\int _0 ^\delta \intxv{\aek \cdot \nabla _v \theta \fek \;\dxv \md t }    \right | \nonumber \\
&\qquad +   \left |\int _0 ^\delta \intxv{a \cdot \nabla _v \theta \;f \;\dxv \md t}    \right | +  \left | \int _\delta ^{+\infty} \!\!\intxv{\left \{ \aek \cdot \nabla _v \theta \;\fek - a \cdot \nabla _v \theta \;f  \right \}\;\dxv \md t}   \right | \nonumber \\
\leq &\, 2 A \delta \|\nabla _v \theta \|_{C^0} \intxv{\fin
\;\dxv} +
\left | \int _\delta ^{+\infty}\!\! \intxv{ (\aek - a) \cdot \nabla _v \theta \;\ind{\{|v|\leq R_0\}}\fek \;\dxv }\md t \right | \nonumber \\
& +  \left | \int _\delta ^{+\infty}\!\! \intxv{ a \cdot \nabla _v
\theta \;(\fek - f)\;\dxv}\md t \right |.
\end{align}
We keep $\delta >0$ fixed and we pass to the limit when $k$ goes
to infinity. Lemma \ref{NonLinTerm} implies that the second term
in the last right hand side can be estimated as
\[
\|\aek - a \|_{\linf (\R^d \times B_{R_0})} = \| a_{\fek} - a_f \|_{\linf (\R^d \times B_{R_0})} \leq C(R_0)
W_1 (\fek (t), f(t)) \to 0 \;\mbox{ when } k \to +\infty
\]
uniformly for $t \in [\delta, T]$, implying, for $T$ large enough
\[
\left | \int _\delta ^{+\infty}\!\! \int _{|v| \leq R_0}{ (\aek -
a) \cdot \nabla _v \theta \fek \;\dxv }\md t \right | \leq
C(R_0)\| \theta \|_{C^1} \int _\delta ^T W_1 (\fek (t), f(t))\;\md
t \to 0
\]
when $k$ goes to infinity. For the third term in the right hand
side of \eqref{EquBil} we use the weak $\star$ convergence $\limk
\fek (t) = f(t)$ in $\mbxv{}$ for any $t\geq \delta$, cf. Proposition \ref{w2properties}
\[
\limk \intxv{a \cdot \nabla _v \theta (\fek (t)- f(t)) \;\dxv } =
0,\;\;t\geq \delta
\]
and we conclude by the Lebesgue dominated convergence theorem
\[
\limk \int _\delta ^{+\infty} \intxv{ a \cdot \nabla _v \theta
(\fek (t,x,v) - f(t,x,v) )\;\dxv}\md t = 0\,.
\]
Passing to the limit in \eqref{EquBil} when $k$ goes to infinity,
we obtain
\[
\limsup _{k \to +\infty} \left |\inttxv{\left \{ \aek \cdot \nabla
_v \theta \fek   -   a \cdot \nabla _v \theta f\right
\}\;\dxv\!\!} \right | \leq 2 A \delta \|\nabla _v \theta \|_{C^0}
\,.
\]
Sending $\delta$ to $0$ we obtain that
\[
\limk \inttxv{ \aek \cdot \nabla _v \theta \; \fek\;\dxv\!\!} =
\inttxv{ a \cdot \nabla _v \theta \; f\;\dxv\!\!}\,.
\]
\end{proof}


\section{Diffusion models}\label{DiffMod}

We intend to introduce a formalism which will allow us to
investigate in a simpler manner the asymptotic behavior of
\eqref{Equ10} and \eqref{Equ31}. This method comes from
gyrokinetic models in plasma physics: when studying the magnetic
confinement we are looking for averaged models with respect to the
fast motion of particles around the magnetic lines. The analysis
relies on the notion of gyro-average operator \cite{BosTraEquSin},
which is a projection onto the space of slow time depending
functions. In other words, projecting means smoothing out the
fluctuations with respect to the fast time variable, corresponding
to the high cyclotronic frequency. This projection appears like a
gyro-average operator. Here the arguments are developed at a
formal level.

We first introduce rigorously the projected measure on the sphere
$r\sphere$ for general measures. Let $f \in \mbxv{}$ be a non
negative bounded measure on $\R^d \times \R^d$. We denote by
$\ave{f}$ the measure corresponding to the linear application
\[
\psi \to \intxv{\psi(x,v)\,\ind{v = 0}  f(x,v)\,\dxv } +
\intxv{\psi\left ( x , r \vsv \right ) \ind{v \neq 0} f(x,v)\,\dxv
}\,,
\]
for all $\psi \in \czcxv$, {\it i.e.,}
\[
\intxv{\psi(x,v) \ave{f}(x,v)\,\dxv} = \int _{v = 0} \psi(x,v)
f(x,v)\,\dxv + \int _{v \neq 0} \psi \left ( x , r \vsv \right )
f(x,v)\,\dxv\,,
\]
for all $\psi \in \czcxv$. Observe that $\ave{f}$ is a non
negative bounded measure,
$$
\intxv{\;\;\ave{f}(x,v)\,\dxv} = \intxv{\;\;f(x,v)\,\dxv},
$$
with $\supp \ave{f} \subset \R^d \times (\A)$. We have the
following characterization.

\begin{pro} \label{VarChar}
Assume that $f$ is a non negative bounded measure on $\R^d \times
\R^d$. Then $\ave{f}$ is the unique measure $F$ satisfying $\supp
F \subset \R^d \times (\A)$,
\[
\int _{v\neq 0} \psi \left ( x , r \vsv \right )F(x,v)\,\dxv =
\int _{v \neq 0}\psi \left ( x , r \vsv \right
)f(x,v)\,\dxv,\;\;\psi \in \czcxv{}
\]
and $F = f$ on $\R^d \times \{0\}$.
\end{pro}
\begin{proof}
The measure $\ave{f}$ defined before satisfies the above
characterization. Indeed, $\supp \ave{f} \subset \R^d \times
(\A)$. Taking now $\psi (x,v) = \varphi (x) \chi (|v|/\delta)$
with $\varphi \in C^0 _c (\R^d)$ and $\delta >0$ one gets
\begin{align*}
\intxv{\varphi (x) \chi \left ( \frac{|v|}{\delta} \right )
\ave{f}(x,v)\,\dxv } = &\,\int _{v = 0} \varphi (x) f(x,v)\,\dxv
\\ &+ \int _{v \neq 0} \varphi (x) \chi \left ( \frac{|v|}{\delta}
\right ) f(x,v)\,\dxv.
\end{align*}
Passing to the limit for $\delta \searrow 0$ yields
\[
\int _{v = 0} \varphi (x) \ave{f}(x,v)\,\dxv = \int _{v = 0}
\varphi (x) f(x,v)\,\dxv,\;\;\varphi \in \czc{}
\]
meaning that $\ave{f} = f$ on $\R^d \times \{0\}$. Therefore one
gets for any $\psi \in \czcxv{}$
\begin{align*}
\int_{v \neq 0} \psi \left ( x , r \vsv \right )\ave{f}(x,v)\,\dxv & =  \int_{|v| = r} \psi  ( x , v )\ave{f}(x,v)\,\dxv \\
& = \int _{v \neq 0} \psi (x,v) \ave{f}(x,v)\,\dxv \\
& = \intxv{\psi \ave{f}}(x,v)\,\dxv- \int _{v = 0}\psi \ave{f}(x,v)\,\dxv \\
& = \intxv{\psi \ave{f}}(x,v)\,\dxv- \int _{v = 0}\psi f(x,v)\,\dxv \\
& = \int _{v \neq 0}  \psi \left ( x , r \vsv \right )
f(x,v)\,\dxv.
\end{align*}
Conversely, let us check that the above characterization exactly defines the measure $\ave{f}$. For any $\psi \in \czcxv{}$ we have
\begin{align*}
\intxv{\psi (x,v) F(x,v)\,\dxv} & = \int _{v = 0} \psi F(x,v)\,\dxv + \int _{v \neq 0} \psi F(x,v)\,\dxv \\
& = \int _{v = 0} \psi (x,v) f(x,v)\,\dxv + \int _{v \neq 0} \psi \left ( x, r \vsv \right ) F(x,v)\,\dxv \\
& = \int _{v = 0} \psi (x,v) f(x,v)\,\dxv + \int _{v \neq 0} \psi
\left ( x, r \vsv \right )f(x,v)\,\dxv
\end{align*}
saying that $F = \ave{f}$.
\end{proof}

By Proposition \ref{VarChar} it is clear that $\ave{\cdot}$ leaves
invariant the measures with support in $\R^d \cup (\A)$. 
Consider $f \in \mbxv{}$. We say that $\Divv \{f \abv \} \in {\cal
M}_b (\R^d \times \R^d)$ if and only if there is a constant $C>0$
such that
\[
\intxv{\abv \cdot \nabla _v \psi f(x,v)\,\dxv } \leq C \|\psi
\|_{\linf{}},\;\;\psi \in \cocxv{}.
\]
In this case there is a bounded measure $\mu$ such that
\[
- \intxv{\abv \cdot \nabla _v \psi f(x,v)\,\dxv } = \intxv{\psi
\mu },\;\;\psi \in \cocxv{}.
\]
By definition we take $\Divv \{f \abv \} = \mu$. The main
motivation for the construction of the projection $\ave{\cdot}$ is
the following result.

\begin{pro}
\label{ZeroAve} For any $f \in \mbxv{}$ such that $ \Divv \{f \abv
\}\in {\cal M}_b (\R^d \times \R^d)$ we have $\ave{\Divv \{f \abv
\}} = 0$.
\end{pro}

\begin{proof}
Let us take $\Divv \{f \abv \} = \mu$. We will check that the zero
measure $0$ satisfies the characterization of $\ave{\mu}$ in
Proposition \ref{VarChar}. Clearly $\supp 0 = \emptyset \subset
\R^d \times (\A)$. For any $\varphi (x) \in \czc{}$ we have
\begin{align*}
\int _{v = 0} \varphi (x) \mu(x,v)\,\dxv & = \limd \intxv{\varphi (x) \chi \left ( \frac{|v|}{\delta} \right ) \mu(x,v)\,\dxv} \\
& = - \limd \intxv{\varphi (x) \chi ^{\;\prime} \left (
\frac{|v|}{\delta} \right )\frac{|v|}{\delta} ( \alpha - \beta
|v|^2) f(x,v)\,\dxv}= 0
\end{align*}
by dominated convergence, since
\[
\left |  \chi ^{\;\prime} \left ( \frac{|v|}{\delta} \right )\frac{|v|}{\delta} ( \alpha - \beta |v|^2)  \right |\leq \alpha \sup _{u \geq 0} |\chi ^{\;\prime} (u) u | + \beta \delta ^2 \sup _{u \geq 0} |\chi ^{\;\prime} (u) u ^3|.
\]
Therefore we deduce that $\Divv \{f \abv\} = 0$ on $\R ^d \times
\{0\}$. Consider now $\psi \in \cocxv{}$ and lets us compute
\begin{align*}
\int _{v \neq 0} \psi \left (x,r \vsv   \right ) & \mu(x,v)\,\dxv = \limd \intxv{\psi \left (x,r \vsv   \right )  \left ( 1 - \chi \left ( \frac{|v|}{\delta} \right ) \right )  \mu(x,v)\,\dxv} \\
& = \limd \intxv{\psi \left (x,r \vsv   \right )  \chi
^{\;\prime}\left ( \frac{|v|}{\delta} \right )  \frac{|v|}{\delta}
(\alpha - \beta |v|^2) f(x,v)\,\dxv} = 0
\end{align*}
since $v \cdot \nabla _v \{ \psi (x, r \vsv)\} = 0$. By density, the same conclusion holds true for any $\psi \in \czcxv{}$ and thus $\ave{\Divv \{f \abv \}} = 0$.
\end{proof}

\begin{remark} \label{SimplerAve}
When $f \in \mbxv{}$ does not charge $\R^d \times \{0\}$,
$\ave{f}$ is given by
\[
\supp \ave{f} \subset \R ^d \times r\sphere,\;\;\int _{v \neq 0}
\psi \left ( x, r \vsv \right ) \ave{f} = \int _{v \neq 0} \psi
\left ( x, r \vsv \right ) f,\;\;\psi \in \czcxv{}
\]
or equivalently
\begin{equation}
\label{Equ34} \intxv{\psi \ave{f} } = \int _{v \neq 0} \psi \left ( x, r \vsv \right ) f,\;\;\psi \in \czcxv{}.
\end{equation}
\end{remark}

\

Using Proposition \ref{ZeroAve} we can obtain, at least formally,
the limit model satisfied by $f = \lime \fe$. By \eqref{Equ2} we
know that $\supp f \subset \R ^d \times (\A)$. The time evolution
of $f$ comes by eliminating $\fo$ in \eqref{Equ3}. For that it is
sufficient to project on the subspace of the measures satisfying
the constraint \eqref{Equ2}, {\it i.e.,} to apply $\ave{\cdot}$.
\begin{equation}
\label{Equ35} \ave{\partial _t f } + \ave{\Divx \{f v\}} + \ave{\Divv \{ f a \}} = 0.
\end{equation}
It is easily seen that $\ave{\partial _t f } = \partial _t \ave{f}
= \partial _t f $ since $\supp f \subset \R ^d \times (\A)$ and
therefore $\ave{f} = f$. We need to compute the last two terms in
\eqref{Equ35}. We show that

\begin{pro}\label{TransportAve}
Assume that $a = a(x)$ is a bounded continuous field. Then we have
the following equalities
\[
\ave{\Divx \{f v \}} = \Divx \{f v \}\;\;\mbox{ if } \;\supp f
\subset \R ^d \times (\A)
\]
\[
\ave{\Divv \{ f a \}} = \Divv \left \{   f \imvv a \right \}
\;\;\mbox{ if } \;\supp f \subset \R ^d \times r\sphere.
\]
As a consequence, \eqref{Equ35} yields the transport equation
\eqref{Equ22} obtained rigorously in Theorems
{\rm\ref{MainResult}} and {\rm\ref{MainResult2}}.
\end{pro}

\begin{proof}
For any $\psi \in \cocxv{}$ we have
\begin{align*}
\intxv{\psi & \ave{\Divx\{f v \}}}  =  \int _{v = 0} \psi \Divx\{f v \} + \int _{v \neq 0} \psixv \Divx\{f v \}\\
& = \limd \intxv{\psi \chivd \Divx\{f v \}} + \limd \intxv{\psixv \left ( 1 - \chivd \right ) \Divx\{f v \}} \\
& = - \limd \intxv{v \cdot \nabla _x \psi \chivd f} - \limd \intxv{v \cdot \nabla _x \psixv \left ( 1 - \chivd \right ) f} \\
& = - \int _{v = 0} v \cdot \nabla _x \psi  f  - \int _{v \neq 0} v \cdot \nabla _x \psixv f  \\
& = - \intxv{v \cdot \nabla _x \psi f } = \intxv{\psi \Divx\{f v
\}}
\end{align*}
saying that $\ave{\Divx \{f v \}} = \Divx \{f v \}$. Assume now
that $\supp f \subset \R ^d \times r\sphere$. It is easily seen
that $\Divv (fa)$ does not charge $\R ^d \times \{0\}$. Indeed, for
any $\psi \in \czcxv{}$ we have by dominated convergence
\begin{align*}
\int _{v = 0} \psi \Divv (fa) & = \limd \intxv{\psi \chivd \Divv (fa)} \\
& = - \limd \intxv{a \cdot \nabla _v \psi \chivd f} - \limd
\intxv{a \cdot \frac{v}{|v|} \frac{1}{\delta} \chipvd \psi f } =
0.
\end{align*}
Therefore we can use \eqref{Equ34}
\begin{align*}
\intxv{\psi \ave{\Divv (fa)}} & = \int _{v \neq 0} \psixv \Divv (fa) \\
& = \limd \intxv{\left ( 1 - \chivd \right ) \psixv \Divv (fa) }\\
& = - \limd \intxv{\left ( 1 - \chivd \right ) \frac{r}{|v|} \imvv a \cdot (\nabla _v \psi ) \left ( x, r\vsv \right ) f } \\
& \quad + \limd \intxv{\;\;\frac{1}{\delta} \chipvd \frac{v}{|v|} \cdot a \psixv f } \\
& = - \int _{v \neq 0} \imvv a \cdot \nabla _v \psi f =
\intxv{\psi \;\Divv \left \{f \imvv a  \right \}}.
\end{align*}
\end{proof}

We investigate now the limit when $\eps \searrow 0$ of the
diffusion model \eqref{Equ31}. We are done if we compute
$\ave{\Delta _v f}$ for a non negative bounded measure with
support contained in $\R^d \times r\sphere$. As before we can
check that $\Delta _v f$ does not charge $\R^d \times \{0\}$ and
therefore, thanks to \eqref{Equ34}, we obtain after some
computations
\begin{equation}
\label{Equ37} \intxv{\psi \ave{\Delta _v f}} = \int _{v \neq 0} \psixv \Delta _v f = \int _{v \neq 0} \Delta _v \left \{ \psixv \right \}f,\;\;\psi \in \ctcxv{}.
\end{equation}

\begin{lemma}
\label{ZeroHom} For any function $\varphi \in C^2 (\R^d \setminus
\{0\})$ and any $r >0$ we have
\[
\Delta _v \left \{ \varphi \left ( r \vsv \right ) \right \} = \left ( \frac{r}{|v|}\right ) ^2 \imvv : \partial ^2 _v \varphi \left ( r \vsv \right ) - 2 \frac{r}{|v|} \frac{v \cdot \nabla _v \varphi \left ( r \vsv \right ) }{|v|^2},\;\;v \neq 0.
\]
\end{lemma}

\

\noindent Combining \eqref{Equ37}, Lemma \ref{ZeroHom} and the
fact that $\supp f \subset \R ^d \times r \sphere$ we obtain
\begin{align*}
\intxv{\psi (x,v) \ave{\Delta _v f} } & = \int _{v \neq 0} \left [ \imvv : \partial _v ^2 \psi (x,v)  - 2 \frac{v \cdot \nabla _v \psi (x,v)}{|v|^2} \right]f \nonumber \\
& = \intxv{\psi (x,v) \Divv \left \{ \Divv \left [ f \imvv \right
] + 2 f \frac{v}{|v|^2}  \right \}}. \nonumber
\end{align*}
We deduce the formula
\[
\ave{\Delta _v f} = \Divv \left \{ \Divv \left [ f \imvv \right ] + 2 f \frac{v}{|v|^2}  \right \}
\]
for any $f$ satisfying $\supp f \subset \R ^d \times r \sphere$
and the limit of the Vicsek model \eqref{Equ31} when $\eps
\searrow 0$ becomes
\begin{equation}\label{equnew}
\partial _t f + \Divx (fv) + \Divv \left \{ f \imvv a \right \} = \Divv \left \{ \Divv \left [ f \imvv \right ] + 2 f \frac{v}{|v|^2}  \right \}
\end{equation}
with the initial condition $f(0) = \ave{\fin}$, as stated in
\eqref{Equ22Diff}.


\appendix

\section{Spherical coordinates and the Laplace-Beltrami operator}
\label{A}

In this appendix, we show the computations to relate the equations
written in original variables $(x,v)$ to the equations in
spherical coordinates $(x,\omega)$. Our limit densities have their
support contained in $\R^d \times r \sphere$ and thus reduce to
measures on $\R^d \times r\sphere$. For example, let us consider
the measure on $\R^d \times r\sphere$ still denoted by $f$, given
by
\[
\intxo{\psi (x, \omega) f(x,\omega)\,\mathrm{d}(x,\omega)} = \int
_{v \neq 0} \psixv{} f(x,v)\,\dxv
\]
for any function $\psi \in \czcxo{}$. In particular, to any $f \in
\mbxv{}$ not charging $\R^d \times \{0\}$ it corresponds $\ave{f}
\in \mbxv{}$, with $\supp \ave{f} \subset \R ^d \times r \sphere$,
whose characterization is
\[
\intxo{\psi (x, \omega)\ave{f}(x,\omega)\,\mathrm{d}(x,\omega)} =
\int _{v \neq 0}\psixv f(x,v)\,\dxv.
\]
We intend to write the previous limit models (in Theorems
\ref{MainResult}, \ref{MainResult2}, and \eqref{equnew}) in
spherical coordinates.

\begin{pro}
\label{SpherCoord} Assume that $f \in \mbxv{}$, $\supp f \subset
\R^d \times r \sphere$ and let us denote by $F \in \mbxo $ its
corresponding measure on $\R^d \times r \sphere$. Therefore we
have
\[
\ave{\Divx (fv)} = \Divx (F \omega),\;\;\ave{\Divv (fa)} = \Divo
\left \{F \imoo a   \right \},\;\;\ave{\Delta _v f } = \Delta
_\omega F.
\]
\end{pro}
\begin{proof}
Thanks to Proposition \ref{TransportAve} we have for any $\psi \in
\cocxo{}$
\begin{align*}
\intxo{\psi (x, \omega) \ave{\Divx (fv)}} & = \int _{v \neq 0} \psixv \Divx (fv) = - \int _{v\neq 0} v \cdot \nabla _x \psixv f\\
& = - \int _{v \neq 0} r \vsv \cdot \nabla _x \psixv f = -
\intxo{\omega \cdot \nabla _x \psi (x, \omega) F}
\end{align*}
and thus $\ave{\Divx (fv)} = \Divx (F \omega)$. Similarly we can
write
\begin{align*}
\intxo{\psi (x, \omega) \ave{\Divv (fa)}} & = \int _{v \neq 0} \psixv \ave{\Divv (fa)}(\mathrm{d}(x,v)) \nonumber \\
& = \int _{v \neq 0} \psixv \Divv \left \{f \imvv a\right \} \nonumber \\
& = - \int _{v \neq 0} \frac{r}{|v|}\imvv a \cdot \imvv \nabla _v\psixv f \nonumber \\
& = - \int _{v \neq 0} \imvv a \cdot \imvv \nabla _v \psixv f \nonumber \\
& = - \intxo{\imoo a \cdot \imoo \nabla _v \psi (x, \omega) F}\nonumber \\
& = - \intxo{\imoo a \cdot \nabla _\omega  \psi (x, \omega) F}\nonumber \\
& = \intxo{\psi (x, \omega) \Divo \left \{ F\imoo a\right
\}}\nonumber
\end{align*}
and therefore
\[
\ave{\Divv (fa)} = \Divo \left \{F \imoo a   \right \}.
\]
Here $\Divo $ stands for the divergence along $r\sphere$ (notice
that $\imoo a$ is a tangent field of $r\sphere$) and $\nabla
_\omega = \imoo \nabla _v$ is the gradient along $r\sphere$. For
the last assertion we appeal to the following well known result
asserting that the Laplace-Beltrami operator coincides with the
Laplacian of the degree zero homogeneous extension, see also
\cite{BCC12}.

\begin{pro}
\label{LaplaceBeltrami} Consider $\varphi = \varphi (\omega)$ a
$C^2$ function on $r\sphere$ and we denote by $\Phi = \Phi (v)$
its degree zero homogeneous extension on $\R^d \setminus \{0\}$
\[
\Phi (v) = \varphi \left ( r \vsv \right ),\;\;v \neq 0.
\]
Therefore we have for any $\omega \in r \sphere$
\[
\Delta _\omega \varphi (\omega) = \Delta _v \Phi (\omega).
\]
\end{pro}

Let us come back to the proof of Proposition \ref{SpherCoord}. For
any $\psi \in C^2_c (\R ^d \times r \sphere)$ we introduce its
degree zero homogeneous extension $\Psi (x,v) = \psixv$. Thanks to
Proposition \ref{LaplaceBeltrami} we can write
\begin{align*}
\intxo{\psi (x,\omega) \ave{\Delta _v f} }& = \int _{v \neq 0}
\psixv  \ave{\Delta _v f} \nonumber  = \int _{v \neq 0} \Psi (x,v)
\Delta _v f = \int _{v \neq 0} \Delta _v \Psi f \\
& = \int _{|v| = r} \Delta _\omega \psi (x,v) f = \intxo{\Delta
_\omega \psi (x,\omega) F} = \intxo{\psi (x,\omega) \Delta _\omega
F}
\end{align*} meaning that $\ave{\Delta _v f } = \Delta _\omega
F$.
\end{proof}

\vskip 6pt

For the sake of completeness, we finally write the equations in
spherical coordinates in $\R^3$. We introduce the spherical
coordinates $\omega = r (\cos \theta \cos \varphi, \cos \theta
\sin \varphi, \sin \theta)$ with the angle variables $(\theta,
\varphi ) \in ]-\pi/2, \pi/2[ \times [0,2\pi [$, and the
orthogonal basis of the tangent space to $r\sphere$
\[
e_\theta = (- \sin \theta \cos \varphi, - \sin \theta \sin
\varphi, \cos \theta),\;\;e_\varphi = (- \cos \theta \sin \varphi,
\cos \theta \cos \varphi, 0)
\]
with $|e_\theta| = 1,\;|e_\varphi| = \cos \theta$. For any smooth
function $u$ on $r\sphere$ we have
\[
\nabla _\omega u = (\nabla _\omega u \cdot e_\theta) e_\theta +
(\nabla _\omega u \cdot e _\varphi ) \frac{e_\varphi}{\cos ^2
\theta} = \frac{1}{r} \partial _\theta u \;e _\theta +
\frac{1}{r\cos ^2 \theta} \partial _\varphi u \;e _\varphi
\]
and for any smooth tangent field $\xi = \xi _\theta e _\theta +
\xi _\varphi e _\varphi $ we have
\[
\Divo \xi = \frac{1}{r} \left \{\frac{1}{\cos \theta} \partial
_\theta (\xi _\theta \cos \theta) + \partial _\varphi \xi _\varphi
\right \}.
\]
The coordinates of the tangent field $\xi := F \imoo a$ are
$
\xi _\theta = \xi \cdot e _\theta = F a_\theta,\;\;\xi _\varphi =
\frac{\xi \cdot e _\varphi }{\cos ^2 \theta } = F a_\varphi
$
and we obtain
\[
\ave{\Divv (fa)} = \Divo \left \{ F \imoo a \right \} =
\frac{1}{r} \left \{ \frac{1}{\cos \theta} \partial _\theta (F
a_\theta \cos \theta) + \partial _\varphi ( F a_\varphi )  \right
\}.
\]
The spherical Laplacian is given by
\begin{align*}
\Delta _\omega F & = \Divo (\nabla _\omega F) = \frac{1}{r} \left \{\frac{1}{\cos \theta} \frac{\partial}{\partial \theta} \left ( \frac{\cos \theta}{r} \partial _\theta F\right ) + \frac{\partial}{\partial \varphi } \left ( \frac{1}{r \cos ^2 \theta} \partial _\varphi F \right )    \right \}\\
& = \frac{1}{r^2}\left \{ \frac{1}{\cos \theta}
\frac{\partial}{\partial \theta} ( \cos \theta \;\partial _\theta
F ) + \frac{1}{\cos ^2 \theta} \;\partial ^2 _\varphi F   \right
\}.
\end{align*}
\begin{pro}
The limit transport equation obtained in \eqref{equnew} for $\R^3$
is
\[
\partial _t F + \omega \cdot \nabla _x F + \frac{1}{r} \left \{\frac{\partial _\theta (F a_\theta \cos \theta )}{\cos \theta}  + \partial _\varphi ( F a _\varphi )   \right \} = \frac{1}{r^2} \left \{\frac{1}{\cos \theta} \frac{\partial}{\partial \theta} ( \cos \theta \;\partial _\theta F ) + \frac{1}{\cos ^2 \theta} \;\partial ^2 _\varphi F    \right \}.
\]
\end{pro}

We recall here the proof of Proposition \ref{LaplaceBeltrami}. It
is a consequence of a more general result.

\begin{pro}
\label{MoreGenRes} Let us consider a function $\varphi = \varphi
(v) \in C^2 (\R^d)$, $d \geq 2$ which writes in polar coordinates
$ \varphi (v) = \tvarphi (\rho, \sigma),\;\;\rho = |v|
>0,\;\;\sigma = \vsv \in \sphere. $ Therefore for any $v \neq 0$ we
have
\[
\Delta _v \varphi (v) = \frac{1}{\rho ^{N-1}} \frac{\partial}{\partial \rho} ( \rho ^{N-1} \partial _\rho \tvarphi ) + \frac{1}{\rho ^2 } \Delta _\sigma \tvarphi (\rho, \sigma),\;\;\rho = |v| >0,\;\;\sigma = \vsv.
\]
\end{pro}
\begin{proof}
Consider a smooth function $\psi = \psi (v) \in C^2$ with compact
support in $\R ^N \setminus \{0\}$, which writes in polar
coordinates $ \psi (v) = \tpsi (\rho, \sigma),\;\;\rho = |v|
>0,\;\;\sigma = \vsv \in \sphere. $ We have
\[
\frac{\partial \tvarphi }{\partial \rho } = \nabla _v \varphi \cdot \sigma,\;\;\nabla _v \varphi = (\nabla _v \varphi \cdot \sigma ) \sigma + (I - \sigma \otimes \sigma) \nabla _v \varphi =
\frac{\partial \tvarphi }{\partial \rho }\;\sigma + \nabla _{\omega = \rho \sigma} \tvarphi
\]
and
\[
\frac{\partial \tpsi }{\partial \rho } = \nabla _v \psi \cdot \sigma,\;\;\nabla _v \psi = (\nabla _v \psi \cdot \sigma ) \sigma + (I - \sigma \otimes \sigma) \nabla _v \psi =
\frac{\partial \tpsi }{\partial \rho }\;\sigma + \nabla _{\omega = \rho \sigma} \tpsi.
\]
Integrating by parts yields
\begin{align*}
- \intvN{\Delta _v \varphi \;\psi (v)} & =  \intvN{\nabla _v
\varphi \cdot \nabla _v \psi } = \int_{\R_+} \int _{S^{N-1}} \left
\{ \frac{\partial \tvarphi}{\partial \rho} \frac{\partial
\tpsi}{\partial \rho} + \frac{1}{\rho ^2} \nabla _\sigma \tvarphi
\cdot \nabla _\sigma \tpsi \right \}
\;\mathrm{d}\sigma \rho ^{N-1} \;\mathrm{d}\rho  \nonumber \\
& = - \int _{S^{N-1}} \int _{\R_+} \tpsi \frac{\partial}{\partial
\rho } \left ( \rho ^{N-1} \frac{\partial \tvarphi}{\partial \rho}
\right )
\;\mathrm{d}\rho\;\mathrm{d}\sigma - \int _{\R_+} \frac{\rho ^{N-1}}{\rho ^2} \int _{S^{N-1}} \tpsi \;\Delta _\sigma \tvarphi \;\mathrm{d}\sigma \;\mathrm{d}\rho \nonumber \\
& = - \intvN{\psi (v) \left \{ \frac{1}{\rho ^{N-1}}
\frac{\partial}{\partial \rho} ( \rho ^{N-1} \partial _\rho
\tvarphi ) + \frac{1}{\rho ^2 } \Delta _\sigma \tvarphi    \right
\}}
\end{align*}
and therefore
\[
\Delta _v \varphi (v) = \frac{1}{\rho ^{N-1}} \frac{\partial}{\partial \rho} ( \rho ^{N-1} \partial _\rho \tvarphi ) + \frac{1}{\rho ^2 } \Delta _\sigma \tvarphi (\rho, \sigma),\;\;\rho = |v| >0,\;\;\sigma = \vsv.
\]
\end{proof}
\begin{proof} (of Proposition \ref{LaplaceBeltrami})
The degree zero homogeneous extension $\Phi (v) = \varphi \left (
r \vsv \right )$ does not depend on the polar radius $\Phi (v) =
\tilde{\Phi} (\sigma) = \varphi (\omega = r \sigma),\;\;\sigma =
\vsv.$ Thanks to Proposition \ref{MoreGenRes}, we deduce
$
\Delta _v \Phi = \frac{1}{\rho ^2} \Delta _\sigma \tilde{\Phi} = \frac{r^2}{\rho ^2} \Delta _\omega \varphi .
$
Taking $\rho = r $, which means $v = r \sigma = \omega$ we obtain
$
\Delta _v \Phi (\omega) = \Delta _\omega \varphi
(\omega),\;\;\omega \in r \sphere.
$
\end{proof}

\subsection*{Acknowledgments}
JAC was supported by projects MTM2011-27739-C04-02 and
2009-SGR-345 from Ag\`encia de Gesti\'o d'Ajuts Universitaris i de
Recerca-Generalitat de Catalunya.

\end{document}